\documentclass[12pt]{amsart}
\usepackage{amssymb}
\usepackage{amsfonts}
\usepackage{amsmath}
\usepackage{array}
\usepackage{rotating}
\usepackage[utf8]{inputenc}
\usepackage{enumitem}
\usepackage{placeins}
\usepackage{pspicture}
\usepackage{epsf}

\usepackage{xcolor}

\usepackage{tikz}

\usepackage[matrix,ps,xdvi]{xypic} 
\newdimen\vcadre\vcadre=0.1cm 
\newdimen\hcadre\hcadre=0.1cm 
\def\GrTeXBox#1{\vbox{\vskip\vcadre\hbox{\hskip\hcadre%
      $#1$%
   \hskip\hcadre}\vskip\vcadre}}
\def\arx#1[#2]{\ifcase#1 \relax \or%
  \ar @{-}[#2]  \or%
  \ar @2{-}[#2] \or%
  \ar @{--}[#2] \or%
  \ar @2{.}[#2] \or%
  \ar @{~}[#2]  \fi}

\newtheorem{example}{Example}[section]

\newtheorem{note}[example]{Note}
\newtheorem{theorem}[example]{Theorem}

\newtheorem{corollary}[example]{Corollary}

\newtheorem{definition}[example]{Definition}
\newtheorem{proposition}[example]{Proposition}
\newtheorem{algorithm}[example]{Algorithm}
\newtheorem{lemma}[example]{Lemma}

\def\Proof{\noindent \it Proof -- \rm}
\def\qed{\hspace{3.5mm} \hfill \vbox{\hrule height 3pt depth 2 pt width 2mm}
\bigskip}

\def\O{{\bf 0}}
\def\I{{\bf 1}}

\def\Hh{{\mathcal H}}
\def\rk{{\rm rk}}
\def\pf{{\rm pf}}

\def\PQSym{{\bf PQSym}}
\def\CQSym{{\bf CQSym}}

\def\ev{{\rm Ev}}

\def\Park{{\rm Park}}

\def\NC{{\rm NC}}

\def\<{\langle}
\def\>{\rangle}

\def\F{{\bf F}}

\def\P{{\bf P}}

\def\SG{{\mathfrak S}}

\def\a{{\bf a}}
\def\b{{\bf b}}

\def\NCSF{{\bf Sym}}
\def\QSym{{\it QSym}}

\def\BT{{\rm BT}}

\usepackage{placeins}
\usepackage{pspicture}

\def\PF{{\rm PF}}

\def\NDPF{{\rm NDPF}}

\def\Tabvrule{\vrule width-0.4pt}       
\def\Tabhrule{\hrule \hrule height-0.4pt} 
\def\Tabstrut{\vrule height2.2ex 
                     depth0.8ex  
                     width0ex    
\relax}

\def\PasCase#1{\omit%
            $\vcenter{\hbox {\vbox to 0.4pt{}}
               \hbox{\makebox[3ex]{\Tabstrut$#1$}}}%
               \Tabvrule$}
\def\PasCasePoint{\PasCase{\cdot}}
\def\DessinCarre#1{%
    \vcenter{\hbox{}\hrule
             \hbox{\vrule\makebox[3ex]{\Tabstrut$#1$}\vrule}\Tabhrule}%
             \Tabvrule}
\def\GenRuban#1{\vcenter{\halign{&$\DessinCarre{##}$\cr#1}}\egroup}

\def\sTabvrule{\vrule width-0.4pt}
\def\sTabhrule{\hrule \hrule height-0.4pt}
\def\sTabstrut{\vrule height1.6ex depth0.6ex width0ex \relax}
\def\sDessinCarre#1{%
    \vcenter{\hbox{}\hrule
             \hbox{\vrule\makebox[2.3ex]%
                  {\sTabstrut$\scriptstyle#1$}\vrule}\sTabhrule}%
             \sTabvrule}
\def\sGenRuban#1{\vcenter{\halign{&$\sDessinCarre{##}$\cr#1}}\egroup}

\def\ruban{%
  \bgroup
  \let\ =\omit
  \let\\=\cr
  \let\x=\times
  \let\.=\PasCasePoint
  \offinterlineskip
  \GenRuban}

\def\sruban{%
  \bgroup
  \let\ =\omit
  \let\x=\times
  \let\\=\cr
  \offinterlineskip
  \sGenRuban}

\newcommand\Matching[2]{%
  \begin{tikzpicture}
    \foreach \x in {1,...,#1}{
       \draw[circle,fill] (\x,0)circle[radius=1mm]node[below]{$\x$};
    }
    \foreach \x/\y in {#2} {
       \draw(\x,0) to[bend left=30] (\y,0);
    }
  \end{tikzpicture}%
}

\newcommand\Matchingsmall[2]{%
  \begin{tikzpicture}
    \foreach \x in {1,...,#1}{
       \draw[circle,fill] (\x/2.5,0)circle[radius=.7mm]node[below]{$\x$};
    }
    \foreach \x/\y in {#2} {
       \draw(\x/2.5,0) to[bend left=80] (\y/2.5,0);
    }
  \end{tikzpicture}%
}

\def\up#1#2{\put(#1,#2){\circle*{0.9}\Line(2,2)}}
\def\down#1#2{\put(#1,#2){\circle*{0.9}\Line(2,-2)}}
\def\hori#1#2{\put(#1,#2){\circle*{0.9}\Line(2,0)}}
\def\pt#1#2{\put(#1,#2){\circle*{0.9}}}

\def\treetoncpncp{\varphi}
\def\gf#1#2{\genfrac{}{}{0pt}{}{#1}{#2}}

\newdimen\vcadre\vcadre=0.2cm 
\newdimen\hcadre\hcadre=0.2cm 

\def\cerp#1#2{\put(#1,#2){\circle*{0.7}}}
\def\cerg#1#2{\put(#1,#2){\circle*{1}}}

\def\arbtgb{\begin{picture}(3,6)\cerg25\cerp13\cerp33\cerp11
 \put(2,5){\Line(-1,-2)}
 \put(2,5){\Line( 1,-2)}
 \put(1,3){\Line( 0,-2)}
\end{picture}}

\def\cerp#1#2{\put(#1,#2){\circle*{0.7}}}
\def\cerg#1#2{\put(#1,#2){\circle*{1}}}

\title[Noncommutative Symmetric Functions and Lagrange Inversion II]%
{Noncommutative Symmetric Functions and Lagrange Inversion II: \\
Noncrossing partitions and the Farahat-Higman algebra}

\author[J.-C.~Novelli and J.-Y.~Thibon]%
{Jean-Christophe Novelli and Jean-Yves Thibon}

\address[] {Laboratoire d'Informatique Gaspard Monge, Universit\'e Gustave
Eiffel, CNRS, ENPC, ESIEE-Paris, \\
5 Boulevard Descartes \\Champs-sur-Marne \\77454 Marne-la-Vall\'ee cedex 2 \\
FRANCE}
\email[Jean-Christophe Novelli]{novelli@univ-mlv.fr}
\email[Jean-Yves Thibon]{jyt@univ-mlv.fr} 
\date{}



\begin{document}

\begin{abstract}
We introduce a new pair of mutually dual bases of noncommutative symmetric
functions and quasi-symmetric functions, and use it to derive generalizations
of several results on the reduced incidence algebra of the lattice of
noncrossing partitions. As a consequence, we obtain a quasi-symmetric
version of the Farahat-Higman algebra.
\end{abstract}

\maketitle
\section{Introduction}

By the \emph{Lagrange series}, we shall mean the (unique) formal power series 
\begin{equation}
g(t)=\sum_{n\ge 0}g_nt^n
\end{equation} 

solving the functional equation

\begin{equation}\label{eq:lag}
g(t) = f(tg(t)) =
  \sum_{ n\ge 1}f_nt^ng(t)^n\ \text{where}\ f(t)=\sum_{n\ge 0}f_nt^n,\ f_0=1.
\end{equation}

Besides its numerous applications in enumerative combinatorics, where the
$f_n$ are specified numbers, the {\it generic} Lagrange series (where the
$f_n$ are indeterminates) is of great interest in algebraic combinatorics.
Specifically, if one interprets the $f_n$ as the homogenous symmetric
functions $f_n=h_n(X)$, the symmetric function $g_n(X)$

\begin{enumerate}[label=(\roman*)]
\item is the Frobenius characteristic of the representation of the symmetric
group $\SG_n$ on the set $\PF_n$ of parking functions of lenght $n$
\cite{Hai1};
\item provides an isomorphism between the reduced incidence Hopf algebra
$\Hh_\NC$ of the family of lattices of noncrossing partitions and the Hopf
algebra $Sym$ of symmetric functions by sending the class $y_n$ of
$[\O_{n+1},\I_{n+1}]$ to $g_n$ \cite{Einz};
\item provides an isomorphism between the Farahat-Higman algebra of symmetric
groups and symmetric functions by identifying\footnote{Actually, Macdonald
uses the equivalent basis $h_\mu^*(X)=g^\mu(-X)$.} the reduced classes $c_\mu$
with the dual basis of $g^\mu:=g_{\mu_1}\cdots g_{\mu_r}$ \cite{Mcd,GJ1,GJ2}.
\end{enumerate}

The Lagrange series has a natural noncommutative version, already apparent in
the original version of Raney's combinatorial proof \cite{Ran}: if in
\eqref{eq:lag} the $f_n$ are interpreted as non-commuting variables, including
$f_0$, $g_n$ becomes the sum of all Łukasiewicz words of length $n+1$: writing
for short $f_{i_1i_2\cdots}$ for $f_{i_1}f_{i_2}\cdots$,
\begin{equation}
 g_0 = f_0,\
 g_1 = f_{10},\
 g_2 = f_{200}+f_{110},\
 g_3 = f_{3000}+f_{2100}+f_{2010}+f_{1200}+f_{1110},
 \ldots
\end{equation}
{\it i.e.}, the Polish codes for plane rooted trees on $n$ vertices (obtained
by reading the arities of the nodes in prefix order). These words also encode
in a natural way various Catalan sets. Setting $f_i=a^ib$, we obtain Dyck
words (with an extra $b$ at the end). Seeing $f_{i_1i_2\cdots i_r}$ as
encoding the nondecreasing word $1^{i_1}2^{i_2}\cdots r^{i_r}$, we obtain a
nondecreasing parking function, which can itself be decoded as a noncrossing
partition, whose blocks are encoded by their minimal elements repeated as many
times as the lengths of the blocks.

\medskip
{\footnotesize
For example, the  word $f_{2100}$ encodes the plane tree $\arbtgb\ $, the
Dyck word $aababb\cdot b$, the nondecreasing parking function $112$ and the
noncrossing partition $13|2$.
}

\medskip
The noncommutative Lagrange series can be interpreted as a noncommutative
symmetric function:  keeping the functional equation \eqref{eq:lag}, we set
$f_n=S_n=S_n(A)$ with $f_0=1$, and obtain
\begin{equation}
\label{g0123}
 g_0 = 1,\ 
 g_1 = S_1,\
 g_2 = S_2+S^{11},\
 g_3 = S_3+2S^{21}+S^{12}+S^{111},\dots
\end{equation}
and we may ask whether there are analogues  for these noncommutative symmetric
functions of Properties (i), (ii), (iii).

Point (i) has been dealt with in~\cite{NTLag}: $g_n(A)$ is the noncommutative
Frobenius characteristic of the natural representation of the $0$-Hecke
algebra $H_n(0)$ on parking functions. Various consequences of this fact,
including a noncommutative $q$-Lagrange formula and generalisations to
$(k,\ell)$-parking functions have been derived there. Other applications have
been given in \cite{NTDup,NTm,JMNT}.

The aim of the present paper is to investigate points (ii) and (iii) in the
noncommutative setting. Introducing the multiplicative basis
$g^I:=g_{i_1}\cdots g_{i_r}$ of $\NCSF$, and computing its coproduct and
antipode, we obtain natural noncommutative versions of these results.

\begin{definition}
The \emph{ordered type} of a noncrossing partition is the composition formed
by the length of its blocks, ordered by increasing values of their minima. Its
\emph{reduced ordered type} is the composition obtained from its ordered type
by subtracting~1 to its components and removing the zeros.

We define the \emph{ordered cycle type} and the \emph{reduced
ordered cycle type} of a permutation similarly. 
\end{definition}

We have then the following interpretation of the coproduct:

\begin{theorem}
\label{th:coprod}
The coefficient $a_{IJ}$ in the coproduct
\begin{equation}
\Delta g_n = \sum_{I,J}a_{IJ}\ g^I\otimes g^J
\end{equation}
is equal to the number of noncrossing partitions $\pi$ of $[n+1]$ of reduced
ordered type~$I$, and whose (right) Kreweras complement $\pi'$ has reduced
ordered type~$J$.
\end{theorem}

This implies a quasi-symmetric refinement of Macdonald's realization of the
graded Farahat-Higman algebra.
For a composition $I=(i_1,\dots,i_r)$, define the canonical permutation
$\sigma_I$ as the permutation of $\SG_{|I|+r}$ whose nontrivial cycles are
\begin{equation}
(12\dots i_1+1)(i_1+2\dots i_1+i_2+2)\dots
(i_1+\cdots+i_{r-1}+r-1\dots i_1+\dots+i_r+r)).
\end{equation}

\begin{corollary}
Let $c_I\in QSym$ be the dual basis of $g^I$. The coefficient $a_{JK}^I$ in
the product
\begin{equation}
c_Jc_K = \sum_{I}a_{JK}^I c_I
\end{equation}
is equal to the number of minimal factorizations $\sigma_I=\alpha\beta$ of the
canonical permutation~$\sigma_I$ of reduced cycle type $I$ with $\alpha$ of
reduced cycle type $J$ and $\beta$ of reduced cycle type~$K$.
\end{corollary}

While this gives back the result of Macdonald by summing over compositions
with the same underlying partitions, the $a_{JK}^I$ only count factorizations
of particular permutations, and this result is rather to be interpreted as
providing a noncommutative version of the reduced incidence algebra of
the lattices of noncrossing partitions.

\medskip
{\footnotesize
For example, $\Delta g_5$ contains the terms $7g^{12}\otimes g^{11}$
and $11g^{21}\otimes g^{11}$, so any 6-cycle in $\SG_n$  has
$18=7+11$ factorizations into permutations of cycle types
$(321^{n-5})$ and $(221^{n-4})$, but the refined coefficients 7 and
11 are only meaningful for the particular 6-cycle $(123456)$.
}

\medskip
Theorem~\ref{th:coprod} is a noncommutative analogue of the main result of
\cite{Einz}, which establishes an isomorphism of Hopf algebras between the
reduced incidence algebra $\Hh_\NC$ of noncrossing partitions and symmetric
functions.  Another result of \cite{Einz} is a combinatorial  description of
the antipode of $\Hh_\NC$. This amounts to computing $g(-X)$ in the basis
$g^\mu$.

Rather than working with the antipode, we shall work with the automorphism
$S_n(A)\mapsto S_n(-A)=(-1)^n\Lambda_n(A)$, and prove the equivalent result

\begin{theorem}\label{th:2}
Define coefficients $a_I$ by
\begin{equation}
g_n(-A) =  \sum_{I\vDash n}(-1)^{\ell(I)}a_Ig^I
\end{equation}
Then, 
\begin{equation}
a_I = \sum_{J\le 2  I}\<M_J,g_{2n}\> =\<E_{2 I},g_{2n}\>
\end{equation}
where $E$ is the so-called essential basis of quasi-symmetric functions.
It is equal to the number of sylvester classes of words of evaluation
$2\bar I$~\cite{NTDup,NTm}
or alternatively, to the number of parking quasi-ribbons of shape
$(2I)^\sim$~\cite{NTLag},
and also to the number of nondecreasing parking functions of type
$2I+1^r:=(2i_1+1,\ldots,2i_r+1)$, which is the same as the number of plane
trees whose arities of internal nodes read in prefix order form the
composition $2I+1^r$.
\end{theorem}

\medskip
{\footnotesize
For example, the term $5g^{21}$ in 
\begin{equation}
g_3(-A) =-g_3+ 5g^{21}+ 3g^{12}-12 g^{111}
\end{equation}
corresponds to the 5 parking quasi-ribbons of shape
$(42)^\sim=12111$ which are
\begin{equation}
1|23|4|5|6,\ 1|22|4|5|6,\ 1|22|3|5|6,\ 1|22|3|4|6,\ 1|22|3|4|5,
\end{equation}
and the term $3g^{12}$ corresponds to the 3 parking quasi-ribbons of shape
$(24)^\sim=11121$
\begin{equation}
1|2|3|45|6,\ 1|2|3|44|6,\ 1|2|3|44|5.
\end{equation}
The term $5g^{21}$ corresponds also to the 5 sylvester classes of evaluation
$24$, which are those of the words
\begin{equation}
112222,\ 211222,\ 221122, \ 222112, \ 222211,
\end{equation}
which can be read by filling the sectors of the plane trees
of skeleton $53$ as in \cite{HNTtrees}.
}

\medskip
The coefficient $\tilde a_I$ in the antipode
\begin{equation}
\tilde\omega(g_n) = \sum_{I\vDash n}(-1)^{\ell(I)}\tilde a_Ig^I
\end{equation}
also has  an explicit, but more complicated interpretation. 
\begin{equation}
\<c_I,\tilde\omega(g_n)\>
 = (-1)^n\sum_{J\vDash n}\<V_I,g^J\>\<M_{\overline{J}},g\>.
\end{equation}
The factor $\<M_{\overline{J}},g\>$ is a number of nondecreasing parking
functions and the $\<V_I,g^J\>$ count parking quasi-ribbons and have all the
same sign. This is therefore a cancel\-lation-free combinatorial formula.

\medskip
{\footnotesize
For example,
\begin{equation}
\tilde\omega(g_3)=-12g^{111} + 4g^{12} + 4g^{21} - g^{3}
\end{equation}
and the contributions to the coefficient of $g^{21}$ are
\begin{equation}
\<V_{21},g_3\>\<M_3,g\> = -3\times 1
\end{equation}
where the factor $-3$ counts the parking quasi-ribbons  $11|2,11|3,12|3$,
and
\begin{equation}
\<V_{21},g^{21}\>\<M_{12},g\> = -1\times 1
\end{equation}
where, dualizing,
$\<V_{21},g^{21}\>=\<V_2\otimes V_1,g_2\otimes g_1\>=-1\times 1$.

Similarly, the contributions to the coefficient of $g^{12}$ are
\begin{equation}
\<V_{12},g_3\>\<M_3,g\> = -2\times 1,
\end{equation}
the $-2$ counts the parking quasi-ribbons $1|22,1|23$, and
\begin{equation}
\<V_{12},g^{12}\>\<M_{21},g\> = -1\times 2,
\end{equation}
where, dualizing,
$\<V_{12},g^{12}\>=\<V_1\otimes V_2,g_1\otimes g_2\>=-1\times 1$.
}

\subsection*{Byproducts}
The proofs of the aforementioned results rely on a couple of elementary
combinatorial properties that we have not been able to find in the
literature and appear to be of independent interest.

First, given a binary tree $T$ and its infix labeling ({\it i.e.}, its corresponding
binary search tree), we decribe a a straightforward algorithm for visiting
cyclically its nodes ({\it i.e}, going from the node labelled $i$ to that
labelled $i+1\mod n$): move one step down the right branch of $i$ (if $i$ is
at its bottom, then go to the top of the branch), then move one step up the
current left branch (again, if $i$ is at its top, go to the bottom of the
branch).  This property is  easily  proved by observing that any such tree
hides two permutations whose product is the standard long cycle, see
Note~\ref{note-infix}.

Second, we prove that a noncrossing partition can be reconstructed from
its ordered type and the ordered type of its right Kreweras complement, and
provide an algorithm doing this.
This property is the key ingredient in the calculation of $\Delta g_n$,
see Theorem~\ref{bicomp-inj}.

\bigskip
This paper is a continuation of \cite{NTLag}, to which the reader is referred
for background and notation.

\section{The Lagrange bases of $\NCSF$ and $QSym$}

\subsection{The Lagrange basis in $\NCSF$}

The Lagrange series in $\NCSF(A)$ is defined by
\begin{equation}
\label{eq:nclag}
g(A) = 1 + \sum_{n\ge 1}S_n(A)g(A)^n
\end{equation}
and we denote by $g_n$ its homogenous component of degree $n$. If $X=(x_i)$ is
a sequence of mutually commuting variables, $g_n(X)$ becomes an ordinary
symmetric function. It is equal to $h_n^*(-X)$ where $*$ is Macdonald's
involution \cite[Ex. 24 p. 36]{Mcd}.

As  mentioned in the introduction, it was shown in~\cite{NTLag} that
\begin{equation}
\label{g2park}
g_n(A) = \sum_{\pi\in\NDPF(n)} S^{\ev(\pi)},
\end{equation}
where $\NDPF$ is the set of nondecreasing parking functions and $\ev(\pi)$ is the
\emph{evaluation} of $\pi$, that is, the ordered sequence of number of
occurrences of $i$ in $\pi$ for $i\ge 1$. Since by convention $S_0=1$, we can
replace $\ev(\pi)$ by the \emph{packed evaluation}, or \emph{type} 
$t(\pi)$ of $\pi$, which is the
composition obtained by removing the zeros in $\ev(\pi)$.

\medskip
{\footnotesize
For example, there are five nondecreasing parking functions: $111$, $112$,
$113$, $122$, and $123$. Forgetting the trailing zeroes, their respective
evaluations are respectively  $3$, $21$, $201$, $12$, and $111$, so that, 
\begin{equation}
\label{g3onS}
g_3 = S_3 + 2S^{21} + S^{12} + S^{111}.
\end{equation}
}

Since $g_n$ begins with a term $S_n$, their products
$g^I=g_{i_1}\dots g_{i_r}$ are triangular on the $S^I$ hence form a basis of
$\NCSF$.
%

Since the coefficient of $S^J$ in $g_n$ is the number of nondecreasing parking
functions of type $J$, or equivalently the number of noncrossing partitions
of ordered type $J$, 
the coefficient of $S^J$ in $g^I$ is the number of nondecreasing parking
functions of type $J$ having breakpoints at the descents of $I$, or the number
of noncrossing partitions of ordered type $J$ finer than the interval
partition of type $I$.

\medskip
{\footnotesize
Ordering compositions in reverse lexicographic order,
{\it e.g.}, $[3, 21, 12, 111]$ for $n=3$,
the matrix of the $g^J$
on the $S^I$ is
\begin{equation}
\label{matG2S}
\begin{pmatrix}
1 & 0 & 0 & 0 \\
2 & 1 & 0 & 0 \\
1 & 0 & 1 & 0 \\
1 & 1 & 1 & 1 \\
\end{pmatrix}
\end{equation}
}
\subsection{A related basis}
This combinatorial description suggests to introduce
another basis 
\begin{equation}
f^I := \sum_{J\geq I} (-1)^{\ell(I)-\ell(J)} g^J,
\end{equation}
where $J\geq I$ means that $J$ is finer than $I$, or that the descents of $I$
are descents of $J$.

The transition matrix from the $f$ to the $S$ is much simpler: the
coefficients are nonnegative integers and each nondecreasing parking function
contributes to the row indexed by its type and to the column indexed by
its breakpoints.

\medskip
{\footnotesize
For example, at $n=3$, the combinatorial description and the matrix from $f$
to $S$ are as follows:
\begin{equation}
\label{matGp2S}
\begin{pmatrix}
111 &     &     &      \\
112 & 113 &     &      \\
    &     & 122 &      \\
    &     &     & 123  \\
\end{pmatrix}
\text{\ \ \ and\ \ \ }
\begin{pmatrix} 
 1  &     &  .  &  .   \\
 1  &  1  &  .  &  .   \\
 .  &  .  &  1  &  .   \\
 .  &  .  &  .  &  1   \\
\end{pmatrix}
\end{equation}
}
This basis will be investigated in a separate paper in
relation to the quasi-symmetric Farahat-Higman algebra.

\subsection{The dual Lagrange basis in $\QSym$}

We denote by $c_I\in QSym$ the dual basis of $(g^I)$.

By definition of the duality between $\NCSF$ and $\QSym$ the transpose of the
matrix in~\eqref{matG2S} is the matrix of the monomial quasi-symmetric
functions $M_I$  in the basis $c_J$.

\subsection{Some other relevant properties}
The expansions of $g_n$ on the bases $S^I$, $\Lambda^I$ and $R_I$ are given in
\cite{NTLag}. It is also proved in this reference that $g$ is invariant under
the involution $S^I\mapsto S^{I^\sim}$, and that $g(-A)$ satisfies the
functional equation

\begin{equation}
\label{eq:gamma}
g(-A)^{-1}= \sum_{n\ge 0}S_n(A)g(-A)^n,
\end{equation}
that is, $g_n(-A)$ is the image of $S_n(A)$ by the antipode of the Hopf
algebra of noncommutative formal diffeomorphisms of \cite{BFK}.

\subsection{The $k$-Lagrange series}

We shall also need the series $g^{(k)}$, defined by the
functional equation
\begin{equation}\label{eq:defgk}
g^{(k)}=\sum_{n\ge 0}S_n \left[g^{(k)}\right]^{kn}.
\end{equation}
It can also be defined as $g^{(k)}=\phi_k(g)$, where $\phi_k$
is the adjoint of the power-sum plethysm operator $\psi^k:\ M_I\mapsto M_{kI}$
on $QSym$.
 
Recall that a \emph{$k$-parking function} is a word over the positive
integers whose nondecreasing rearrangement $a_1a_2\cdots a_n$ satisfies
$a_i\le k(i-1)+1$. Its $k$-\emph{evaluation} 
$\ev_k(\a)$ is essentially the classical
evaluation of a word but we will here define it as the number
of occurrences of all letters from $1$ to $kn+1$. In particular, the
$k$-evaluation of a nonempty $k$-parking function ends with a sequence of at
least $k$ zeros.

Indeed, we have proved in~\cite{NTLag} that the solution of \eqref{eq:defgk}
where $S_0$ is an indeterminate is
\begin{equation}
g^{(k)}=\sum_{\pi\in\NDPF^{(k)}}S^{\ev_k(\pi)}.
\end{equation}
In particular, if one sends $S_0$ to $1$, the coefficient of $S^I$ in
$g_n^{(k)}$ of degree $n$ is the number of nondecreasing $k$-parking functions
of type $I$.

\medskip
{\footnotesize
For example, setting $h=\phi_2(g)=g^{(2)}$,
\begin{equation}
\label{g2}
h_0+h_1+h_2+h_3+\dots
  = S_0+S_1(h_0+h_1+h_2+\dots)^2+ g_2(h_0+h_1+\dots)^4+g_3(h_0+\dots)^6
\end{equation}
yields, by iterated substitutions
\begin{equation}
h_0 = 1,\ \ \
h_1 = S_1,\ \ \ 
h_2 = S^2 + 2 S^{11},\ \ \
h_3 = S_3 + 4 S^{21} + 2 S^{12} + 5 S^{111},
\end{equation}

The $2$-parking functions of size $3$ are
\begin{equation}
111,\ 112,\ 113,\ 114,\ 115,\ 122,\ 133,\ 123,\ 124,\ 125,\ 134,\ 135,
\end{equation}
and one can check that their types indeed encode the
expansion of $h_3$.
}

\medskip
The $k$-evaluations of $k$-parking functions are generalized Łukasiewicz words.
They are the words $w_1\dots w_{kn+1}$ of length $kn+1$ whose partial sums
$k(w_1+\dots+w_{i})-i$ are always nonnegative except when $i=kn+1$ where the
sum becomes strictly negative. This property is easily translated in terms of
generalized Dyck paths: send $w_i$ to $w_i$ times the step $(1,k)$ followed by
a step $(1,-1)$. The conditions on the evaluations mean that the path stays
weakly above the axis on all steps but the last.

\subsection{Connecting the $k$-Lagrange series}

\medskip
There is a simple but useful connection between $k$-parking functions and
$(k-1)$-parking functions.
Let $\a$ be a $k$-parking function
and let $\a'$ be its largest prefix that is a $(k-1)$-parking function.

In terms of evaluations, this means that $\ev(\a')$ is a prefix of $\ev(\a)$.
If $\a'$ is considered as a $(k-1)$-parking function, the corresponding path
ends at height $-1$, and ends at height $k-1$ as a $k$-parking function. Now,
since each downstep decrements the height by one, one can cut the remainder of
the path of $\a$ the first time it reaches each height from $k-2$ down to $0$.
One then gets a total of $k$ (possibly one-downstep) paths, all encoding a
$k$-parking function.  Conversely, given a $(k-1)$-parking function of length
$i$ and a list of $i$ $k$-parking functions, one obtains a $k$-parking
function by concatenating their evaluations.

\medskip
{\footnotesize
For example, consider the $2$-parking function
$1\,2\,2\,4\,6\,9\,11\,14\,17\,17$.
Its evaluation (up to $21$) is
\begin{equation}
120101001010010020000
\end{equation}
Its largest  $1$-parking prefix is $1224$. It is of size $4$ and its
evaluation (as a  $1$-parking function) is $12010$.
We then remove the prefix and cut the remainder into four parts as
\begin{equation}
1001010010020000 = 100.10100.100.20000,
\end{equation}
that all are evaluations of $2$-parking functions.
}

\medskip
Thus, a $k$-parking function can be uniquely decomposed as
$\a=\a'\b_1\cdots \b_i$ where $\a'$ is its maximal $(k-1)$-parking
prefix of length $i$ and the $\b_j$ are $k$-parking functions.
This translates into the following functional equation:

\begin{lemma}
\label{lem:itergk}
The series $g^{(k)}:=\phi_k(g)$ satisfies 
\begin{equation}
g^{(k)}=\sum_{n\ge 0}g^{(k-1)}_n \left[g^{(k)}\right]^n.
\end{equation}
\end{lemma}
\qed

\medskip
{\footnotesize
For example, setting $h=\phi_2(g)=g^{(2)}$,
\begin{equation}
h_0+h_1+h_2+h_3+\cdots
  = g_0+g_1(h_0+h_1+h_2+\cdots)+ g_2(h_0+h_1+\cdots)^2+g_3(h_0+\cdots)^3
\end{equation}
yields
\begin{equation}
\begin{split}
h_0 = 1,\ \ \
h_1 = g_1 = S_1,\ \ \
h_2 = g_1h_1 + g_2 = S^2 + 2S^{11},\ \ \\
h_3 = g_1h_2 + 2g_2h_1 + g_3 =  S_3 + 4S^{21} + 2S^{12} + 5S^{111},
\end{split}
\end{equation}
where $g_i$ is replaced by its expansion on the $S^I$ as in~\eqref{g0123}.
}

\subsection{Base change from $S$ to $g$}

To compute the change of basis from $S$ to $g$, we proceed as in~\cite{JMNT}. 
In this reference, ``noncommutative free cumulants'' $K_n$ are defined by the
functional equation
\begin{equation}
\sigma_1 = \sum_{n\ge 0}K_n \sigma_1^n
\end{equation}
and it is proved that
\begin{equation}
K(A) = \sum_{n\ge 0}K_n(A)=g(-A)^{-1}.
\end{equation}
Setting $g(A)=\sigma_1(B)$, we see that $K_n(B)=S_n(A)$ and that
\begin{equation}
K_n = \sum_{I}k_IS^I \Leftrightarrow S_n =\sum_{I}k_Ig^I.
\end{equation}
To expand $S_n$ on the basis $g^I$, we can therefore apply the recipe given
in~\cite[Eq. (50)]{JMNT}: start from the expansion of $g_{n-1}$ on the
elementary basis, as given in \cite{NTLag}, and replace each $\Lambda^I$ by
$g^{i_1+1,i_2,\ldots,i_r}-g^{1I}$.

\medskip
{\footnotesize
For example, starting with
\begin{equation}
g_3 = \Lambda^3-3\Lambda^{21}-2\Lambda^{12}+5\Lambda^{1111},
\end{equation}
this substitution yields
\begin{equation}
S_4 = (g^4-g^{13})-3(g^{31}- g^{121})-2(g^{22}-g^{112})+5(g^{211}- g^{1111}). 
\end{equation}
The first values are
\begin{align*}
S_1 &=g_1\\
S_2 &= g_2-g^{11}\\
S_3 &= g_3-2g^{21}-g^{12}+2g^{111}\\
S_4 &= g_4-3g^{31}- g^{13}-2 g^{22} +5g^{211}+3 g^{121}+2 g^{112}-5 g^{1111} 
\end{align*}

and one checks that the inverse matrix of~\eqref{matG2S} is indeed
\begin{equation}
\label{matS2G}
\begin{pmatrix}
1  & 0 & 0 & 0 \\
-2 & 1 & 0 & 0 \\
-1 & 0 & 1 & 0 \\
2  &-1 &-1 & 1 \\
\end{pmatrix}
\end{equation}
}

\subsection{An involution}

It is not immediate that there is an analogue of Macdonald's star involution
in the noncommutative setting. Indeed, $g$ does not commute with the $S_n$,
and writing \eqref{eq:nclag} in the (ambiguous) form $g=\sigma_g(A)$ does not
allow to conclude that $g^{-1}=\lambda_{-g}(A)$. However, this relation does
hold, and we have:

\begin{proposition}
The algebra automorphism $F\mapsto \tilde F$ defined on the elementary
symmetric functions by 
\begin{equation}
\Lambda_n\mapsto \tilde\Lambda_n = g_n
\end{equation}
is an involution of $\NCSF$.
\end{proposition}

\Proof
The noncommutative free cumulants being given by
\begin{equation}
K_n(A)=\bar g_n(-A) \quad\text{where\ \ $\bar g(A):=g(A)^{-1}$},
\end{equation}
we have therefore
\begin{equation}
\sigma_1(-A) = \sum_{n\ge 0}\bar g_n(A)\sigma_1(-A)^n
\end{equation}
so that
\begin{equation}
(-1)^n\Lambda_n(A)
 = S_n(-A)
 = g_n|_{S_k\mapsto \bar g_k}
 = \sum_{I\vDash n}c_I\bar g^I\
   \text{where}\ g_n=:\sum_{I\vDash n}c_IS^I,
\end{equation}
and if on the one hand we define coefficients $b_J^I$ by 
\begin{equation}
(-1)^{|I|}\Lambda^I=\sum_Jb^I_JS^J, 
\end{equation}
then $\bar g^I=\sum_Jb_J^Ig^J$, and 
\begin{equation}
\Lambda_n = \sum_{I,J} c_Ib_J^Ig^J.
\end{equation}
But on the other hand
\begin{equation}
(-1)^n\sum_{I,J} c_Ib_J^I\Lambda^J = (-1)^n\sum_Ic_I(-1)^{|I|}S^I=g_n.
\qed
\end{equation}

\goodbreak
\subsection{Proof of Theorem \ref{th:2}}

With this at hand, we can compute the first values of $g_n(-A)$:

{\footnotesize
\begin{align*}
g_1(-A) &= -g_1  \\
g_2(-A) &= -g_2 + 3g^{11}  \\
g_3(-A) &= -g_3 + 5g^{21} + 3g^{12} - 12 g^{111}    \\
g_4(-A) &= -g_4 + 7g^{31} + 5g^{22} + 3g^{13} - 25 g^{211}
           - 18 g^{121} - 12 g^{112} + 55g^{1111}.
\end{align*}
}

One can observe that the sum of the absolute values of the coefficients build
up the sequence 1,4,21,126,... \cite[A003168]{Slo} and that the coefficients
refine the triangle \cite[A102537]{Slo}, which occurs in \cite[Sec. 5.3]{NTm}.
This suggests that the coefficient of $\pm g^I$ should count sylvester classes
of packed words of evaluation $2\bar I$.

We propose to show
\begin{equation}
g_n(-A)
 = \sum_{I\vDash n}
   (-1)^{\ell(I)}\left(\sum_{J\le 2 I}\<M_J,g\>\right)g^I
 = \sum_{I\vDash n}
   (-1)^{\ell(I)}\left(\sum_{J\le  I}\<M_J,\phi_2(g)\>\right)g^I\,,
\end{equation}
where $\phi_2$ is the adjoint of $\psi^2:\ M_I\mapsto M_{2I}$ \cite{NTm}.
Equivalently, we want to prove that
\begin{equation}
g_n(-A)=\sum_{I\vDash n}\<M_I,\phi_2(g(-A))\>g^I.
\end{equation}
We start from the expansion
\begin{equation}
g_n(-A) =(-1)^n\sum_{I\vDash n}\<M_I,g\>\Lambda^I.
\end{equation}
Let $V_I$ be the dual basis of $\Lambda^I$.
According to the previous considerations, we can write
\begin{align}
g_n(-A) &= (-1)^n\sum_{I\vDash n}
           \<M_I,g\>
           \sum_{I_1\vDash i_1,\ldots, I_r\vDash i_r}
             \<V_{I_1},g\>\cdots \<V_{I_r},g\>g^{I_1I_2\cdots I_r}\nonumber\\
        &= (-1)^n\sum_{J\vDash n}
            g^J\sum_{I\le J}
               \<M_I,g\>\<\Delta^rV_J,g_{i_1}\otimes\cdots\otimes g_{i_r}\>\\
        &=  (-1)^n\sum_{J\vDash n}
            \left(\sum_{I\le J}\<M_I,g\>\<V_J,g^I\>\right)g^J\nonumber.
\end{align}
We are thus reduced to show
\begin{equation}
\<M_I,\phi_2(g(-A))\> = (-1)^n\sum_{J\le I}\<V_I,g^J\>\<M_J,g\>.
\end{equation}
Summing the right-hand sides multiplied by $S^I$ yields
\begin{align*}
\sum_I(-1)^{|I|}\sum_{J\le I}
  \<V_I,g^J\>\<M_J,g\>S^I
&=\sum_I\sum_{J\le I}\<V_I,g^J\>\<M_J,g\>\Lambda^I(-A)\\
&=\sum_J\<M_J,g\>\sum_{I\ge J}\<V_I,g^J\>\Lambda^I(-A)\\
&=\sum_J\<M_J,g\>g^J(-A).
\end{align*}
Doing the same with the left-hand sides, we have finally to show that 
\begin{equation}
\phi_2(g(-A))=\sum_J\<M_J,g\>g^J(-A),
\end{equation}
or equivalently, that
\begin{equation}
h := \sum_J\<M_J,g\>g^J
\end{equation}
satisfies
\begin{equation}
h = \sum_{n\ge 0}g_nh^n.
\end{equation}
which follows from Lemma \ref{lem:itergk}, since $h$ defined as above is
obtained by substituting $g_n$ to $S_n$ in $g$.

This concludes the proof of Theorem \ref{th:2}.

\subsection{Another argument}
\label{sec:otherarg}

Instead of Lemma \ref{lem:itergk}, we can rely upon the tilde involution.
This leads to a different combinatorial interpretation of the coefficients.

Recall that
\begin{equation}
\begin{split}
g(-A)^{-1}
= \sum_{n\ge 0} S_n(A)g(-A)^n
                \Leftrightarrow g(-A)&=1-\sum_{n\ge 1}S_n(A)g(-A)^{n+1}\\
&:=\sum_{n\ge 0}S_n(B)g(-A)^n,
\end{split}
\end{equation}
setting
$S_1(B)=0$ and $S_n(B)=-S_{n-1}(A)$ for $n\ge 2$. 
Hence, the coefficient of $S^I$ in $g(-A)$ is equal to
\begin{equation}
\<M_I,g(-A)\>=(-1)^{\ell(I)}\<M_{I+1^r},g(A)\>,
\end{equation}
where $I+1^r=(i_1+1,\ldots,i_r+1)$.
Applying the involution $\tilde\Lambda_n=g_n$, and setting $h=g(-A)$, we have
\begin{equation}
\tilde h 
 = \sum_{n\ge 0}\widetilde{S_n(-A)}\tilde h^n 
 = \sum_{n\ge 0}(-1)^n\tilde\Lambda_n \tilde h^n
 = \sum_{n\ge 0}g_n(-\tilde h)^n.
\end{equation}
This is, up to signs, the functional equation for $g^{(2)}=\phi_2(g)$, so that
\begin{equation}
g_n(-A)=(-1)^n\tilde g^{(2)}_n.
\end{equation}
Hence, the coefficient of $g^I$ in $g_n(-A)$ is
\begin{equation}
\<c_I,g_n(-A)\> =(-1)^n\<V_I,g^{(2)}_n\> =(-1)^n\<V_{2I},g_{2n}\>
\end{equation}
for which a combinatorial interpretation in terms of parking quasi-ribbons is
given in
\cite{NTLag}:
\begin{equation}
g_n(A)=\sum_{I\vDash n}(-1)^{n-\ell(I)}c_{I^\sim}\Lambda^I,
\end{equation}
where $c_I$ is the number of parking quasi-ribbons of shape $I$. 

We can also give a third combinatorial interpretation of $g(-A)$.
The dual basis of~$\Lambda^I$ is 
\begin{equation}
V_I=(-1)^{n-\ell(I)}\sum_{J\le I}M_J
\end{equation}
so that the coeffficient of $g^I$ in $g(-A)$ is equal to $\<V_{2I},g\>$,
hence, replacing $A$ by $-A$, to the coefficient of $(-1)^{2n}S^{2I}=S^{2I}$
in $g(-A)$.

We have seen that the coefficient $\delta_I$ of $S^I$ in $g$
and the  coefficient $\lambda_I$ of $S^I$ in $g(-A)$ are related by
\begin{equation}
\lambda_I =(-1)^{\ell(I)}\delta_{i_1+1,i_2+1,\ldots,i_p+1}.
\end{equation}
We have therefore for the absolute value of the coefficient of $g^I$ in
$g(-A)$
\begin{equation}
\sum_{J\le 2I}\<M_J,g\>=\<M_{2I+1^r},g\>
\end{equation}
which is the number of nondecreasing parking functions
of type $(2i_1+1,\ldots,2i_r+1)$, or equivalently, to the number of
plane trees whose arities of the internal nodes read in infix order yield this
composition.

\subsection{The antipode of $g$}

The antipode $\tilde\omega(g)$ can be obtained by a slight adaptation of the
argument of Section \ref{sec:otherarg}.

Let $h=\tilde\omega(g)=\overline{g(-A)}$. Then,
\begin{equation}
\tilde h=\sum_{n\ge 0}\tilde h^n(-1)^ng_n(A)
\end{equation}
which is, up to signs,
\begin{equation}
f = \sum_{n\ge 0}f^n g_n
\end{equation}
whose solution is $f=\chi(g^{(2)})$, where $\chi$ is the involution
$g^I\mapsto g^{\bar I}$.
Hence,
\begin{equation}
\tilde\omega(g_n)=(-1)^n\widetilde{\chi(g_n^{(2)})}.
\end{equation}
The coefficient de $g^I$ in $\tilde\omega(g_n)$ is therefore 
\begin{equation}
\<c_I,\tilde\omega(g_n)\>
 = (-1)^n\sum_{J\vDash n}\<V_I,g^J\>\<M_{\overline{J}},g\>.
\end{equation}
The factor $\<M_{\overline{J}},g\>$ is a number of nondecreasing parking
functions, and $\<V_I,g^J\>$ counts parking quasi-ribbons with a common sign.
This is therefore a cancellation-free combinatorial formula.

\goodbreak
\medskip
{\footnotesize
To compute the antipode of $g_n$ :
\begin{itemize}
\item Express $g_n^{(2)}$ on the basis $g^I$
\begin{equation}
g_3^{(2)}=g_3+2g^{21}+g^{12}+g^{111}
\end{equation}
\item apply the involution $\chi: g^I\mapsto g^{\bar I}$
  and multiply by $(-1)^n$
\begin{equation}
(-1)^3\chi(g_3^{(2)})=-(g_3+2g^{12}+g^{21}+g^{111})
\end{equation}
\item then expand it on the basis $\Lambda^I$
\begin{equation}
(-1)^3\chi(g_3^{(2)})=-(\Lambda^3-4\Lambda^{21}-4\Lambda^{12}+12\Lambda^{111})
\end{equation}
\item and finish by applying the tilde involution $\Lambda^I\mapsto g^I$
\begin{equation}
\tilde\omega(g_3)=-(g_3-4g^{21}-4g^{12}+12g^{111}).
\end{equation}
\end{itemize}
}

\section{Expansion of the coproduct of $g_n$ on the basis $g^I$}

\subsection{Background: the Hopf algebra of nondecreasing parking functions}

One can also rewrite ~\eqref{g2park} in $\NCSF$ as
\begin{equation}
g_n = \sum_{I} \delta_I S^I,
\end{equation}
where $\delta_I$ is the number of nondecreasing parking functions of type $I$.  

\medskip
{\footnotesize
For example,
\begin{equation}
\label{g3onSbis}
g_3 = S_3 + 2S^{21} + S^{12} + S^{111}
\end{equation}
is obtained from $111$, $112$, $113$, $122$, $123$.
}

\medskip
We have defined in~\cite{NTpark} an algebra $\PQSym$ based on symbols $\F_\a$,
where $\a$ runs over all parking functions.
One can show that $\PQSym$ has a Hopf subalgebra $\CQSym$ whose basis is
defined by
\begin{equation}
\P^\pi = \sum_{\a\uparrow =\pi}\F_\a,
\end{equation}
where $\pi$ is any nondecreasing parking function and the sum runs over all
parking functions with the same nondecreasing rearrangement $\pi$.

If one denotes by $t(\pi)$ the packed evaluation of $\pi$, which coincides
with the ordered type of the noncrossing partition encoded by $\pi$, then, the
map $\phi:\ \P^\pi\mapsto S^{t(\pi)}$ is an epimorphism of Hopf
algebras~\cite{NTpark}, and
\begin{equation}
g = \phi(G),\ \text{where}\ \
G := \sum_{\a\in\PF}\F_\a 
  = \sum_{\pi\in\NDPF}\P^\pi
\end{equation}
is the formal sum of all parking functions.

\medskip
{\footnotesize
For example,
\begin{equation}
G_3 = \P^{111} + \P^{112} + \P^{113} + \P^{122} + \P^{123},
\end{equation}
so that one recovers \eqref{g3onS} and~\eqref{g3onSbis} by sending $\P^\pi$ to
$S^{t(\pi)}$.
}

\medskip
Thus, $\Delta g = (\phi\otimes\phi)(\Delta G)$
and one can get $\Delta g$ from $\Delta G$ which is simpler,
since, as we shall see shortly, it has an intermediate
multiplicity-free expression.

\subsection{Computation of the coproduct in $\CQSym$}

The coproduct in $\CQSym$ in the $\P$ basis is given by
\begin{equation}
\Delta \P^\pi
  = \sum_{\gf{\pi=(uv)\uparrow}{u,v\ {\rm nondecreasing}}}
    \P^{\Park(u)}\otimes \P^{\Park(v)}
\end{equation}
where $\Park$ denotes the operation of \emph{parkization} as described in
\cite{NTpark}, and the sum runs over all nondecreasing words $u,v$ such that
the nondecreasing rearrangement of $uv$ is $\pi$.

\medskip
{\footnotesize
For example,
\begin{equation}
\begin{split}
\Delta\P^{1124} &= 1\otimes\P^{1124} +
        \P^1\otimes \left(\P^{112}+\P^{113}+\P^{123}\right) +
        \P^{11}\otimes\P^{12} \\
       & + \P^{12}\otimes\left(\P^{11}+2\P^{12}\right)
         + \left(\P^{112}+\P^{113}+\P^{123}\right)\otimes\P^1 +
        \P^{1124}\otimes1\,.
\end{split}
\end{equation}
}

Now, as an intermediate computation, we could ``forget'' to parkize $u$ and
$v$ and write the coproduct of $\P^\pi$ as the sum of all terms
$\P^u\otimes\P^v$, over all pairs on nondecreasing words such that $uv=\pi$.
This amounts to making the convention $\P^w=\P^{\Park(w)}$ for an arbitrary
nondecreasing word~$w$.

\medskip
{\footnotesize
For example, with this convention, the coproduct $\Delta\P^{1124}$ becomes
\begin{equation}
\begin{split}
\Delta\P^{1124}
     & = 1\otimes\P^{1124} +
         \P^1\otimes\P^{124} + \P^2\otimes\P^{114} + \P^4\otimes\P^{112} \\
     & + \P^{11}\otimes\P^{24}  + \P^{12}\otimes\P^{14}+
         \P^{14}\otimes\P^{12} + \P^{24}\otimes\P^{11} \\
     & + \P^{112}\otimes\P^4 + \P^{114}\otimes\P^2 + \P^{124}\otimes \P^1
         + \P^{1124}\otimes1\,.
\end{split}
\end{equation}
}

\medskip
\begin{note}{\rm
With this convention, if one forgets to parkize all terms, this expression of
$\Delta G_n$ becomes multiplicity-free, since a term $\P^u\otimes\P^v$ can
only come from a $\Delta\P^\pi$ where $\pi$ is obtained by sorting $u\cdot v$.

In other words, $\Delta G_n$ is the sum of terms $\P^u\otimes\P^v$, over
all pairs on nondecreasing words such that $uv$ is a parking function.
}
\end{note}

Define $G^I=G_{i_1}\cdots G_{i_r}$. We shall prove that $\Delta G$ is actually
a sum of terms $G^I\otimes G^J$.

\subsection{Profiles of nondecreasing words}


Any nondecreasing word $w$ admits a minimal factorization into shifted parking
functions
\begin{equation}
w = w_1w_2\cdots w_k
\end{equation}
\emph{i.e.},  each $w_i$ is obtained by shifting a parking function $\a_i$ by
some integer $b_i$, which we write as $w_i =(\a_i)_{b_i}$ and each $w_i$ is of
maximal length.


\medskip
{\footnotesize
For example,
\begin{equation}
\begin{split}
w &= 2336799
   = (1225688)_1
   = (122)_1 \cdot 6799 \\
  &= (122)_1 \cdot (1244)_5
   = (122)_1 \cdot (12)_5 \cdot 99 \\
  &= (122)_1\cdot (12)_5\cdot (11)_8,
\end{split}
\end{equation}
so that $w$ decomposes as
\begin{equation}
2336799 = 233 \cdot 67 \cdot 99
\end{equation}
and the $a_i$s and the $b_i$s can be read above.
}

\begin{definition}
The \emph{profile} $\pf(w)$ of a word $w$ is the pair 
$\binom{s}{c}=\binom{s_1s_2\cdots s_k}{c_1c_2\cdots c_k}$,
%
where $s_i$ is the first letter of $w_i$, that is, $1+b_i$ and $c_i$ its
length.

We shall say a biword is a profile if there is a word $w$ whose profile is
that biword.
\end{definition}

{\footnotesize
On our example, $\pf(w)=\binom{2\, 6\, 9}{3\,2\,2}$.
}

\medskip
There is a simple characterization of profiles:
\begin{lemma}
\label{carac-profils}
A biword $\binom{s_1\dots s_k}{c_1\dots c_k}$ is a profile iff
$s_{i+1}>s_i+c_i$ for all $i\in[1,k-1]$.
\end{lemma}

\Proof
Let $w$ be a nondecreasing word. Decompose it as $w_1\dots w_k$ as above,
and let $c_i$ be the length of $w_i$.

Since $w_i$ and $w_{i+1}$ are different blocks of the decomposition of $w$,
after shifting the suffix of $w$ starting with $w_i$, its largest prefix which
is a parking function will be exactly $w_i$. So the first letter of $w_{i+1}$
has to be far enough from the first letter of $w_i$, more precisely,
$s_{i+1}-s_{i}$ has to be strictly greater than their distance in the word
which is $c_i$, whence the condition.

Conversely, given a biword satisfying the required conditions, it is easy to
exhibit a word with that profile:
\begin{equation}
w = s_1^{c_1} s_2^{c_2} \dots s_k^{c_k}.
\end{equation}
\qed

\subsection{Biprofiles of pairs of nondecreasing words}

Given two nondecreasing words $u$ and $v$, we define their \emph{biprofile}
as the pair $(\pf(u),\pf(v))$.

\begin{lemma}
\label{parking-biprofil}
Let $(u,v)$ be a pair of nondecreasing words of respective profiles
$S=\binom{s_1\dots s_k}{c_1\dots c_k}$ and
$T=\binom{t_1\dots t_\ell}{d_1\dots d_\ell}$.

Rearrange the biword
$\binom{s_1\dots s_k t_1\dots t_\ell}{c_1\dots c_k d_1\dots d_\ell}$ as a
\emph{joint profile}, so that the top line is weakly increasing. If some
$s_i=t_j$, put the biletter of $s_i$ to the left of the one of $t_j$
and write the result as
$\binom{x_1\dots x_{k+\ell}}{y_1\dots y_{k+\ell}}$.

Then, the concatenation $uv$ is a parking function iff
\begin{equation}
\label{eq-parking-bipro}
\forall m\in[1,k+\ell],\ x_m \leq y_1+\dots+y_{m-1}+1.
\end{equation}
In that case, we say that the biprofile is a \emph{parking biprofile}.
\end{lemma}

Note that in general, the joint profile is not a profile.

\medskip
{\footnotesize
For example, let $u=2336799$ and $v=11$.
Then the concatenation of their profiles gets reordered as
\begin{equation}
\binom{2\, 6\, 9\, 1}{3\, 2\, 2\, 2} = \binom{1\, 2\, 6\, 9}{2\, 3\, 2\, 2}.
\end{equation}
We have the inequalities
\begin{equation}
x_1=1\leq 1,\ \
x_2=2\leq 3,\ \
x_3=6\leq 6,\ \
x_4=9 > 8,
\end{equation}
so that $uv$ is not a parking function and indeed, there are only $7$ values
smaller than or equal to $8$ in $uv$.

One can also check that if $u$ is the same and $v=116=11\cdot 6$ then 
the joint profile is $\binom{12669}{23212}$, all inequalities are satisfied
and $uv$ is indeed a parking function.
}

\bigskip
\Proof
Let us write as before $u=(u_1)\cdots (u_k)$ and $v=(v_1)\cdots (v_\ell)$.
Then rearrange $uv$ as blocks  matching the rearranged concatenated biword in
the statement:
\begin{equation}
\label{blocs-w}
w = (a_1)(a_2)\dots (a_{k+\ell})
\end{equation}
where the $(a_i)$ run over all factors of both $u$ and $v$ (the $i$-th
block $a_i$ comes from $u$ if $x_i$ is some $s_i$).
Writing the blocks $a_i$ as words, \eqref{blocs-w} becomes
\begin{equation}
w = (w_{z_0+1}\dots w_{z_1}) \dot (w_{z_1+1}\dots w_{z_2}) \dots
    (w_{z_{k+\ell-1}+1} \dots w_{z_{k+\ell}}),
\end{equation}
where $z_i=y_1+\dots+y_i$ with the convention $z_0=0$. 

Let us now assume that some $x_m>z_{m-1}+1$.
In this case, $w_{z_{m-1}+1}=x_m>z_{m-1}+1$ and it has no letter strictly
smaller on its right since the $x_i$ are weakly increasing and among each
block, the values are weakly increasing too. So $w$ cannot be a parking
function: it has less than $z_{m-1}$ values smaller than or equal to
$z_{m-1}+1$.

Conversely, assume that all $x_m\leq z_{m-1}+1$.
In that case, each letter beginning a block satisfies
$w_{z_{m-1}+1}\leq z_{m-1}+1$.
Now, any $w_{z_i+j}$ with $j\leq y_i$ is at most $z_i+j$ since the subword
$w_{z_i+1}\dots w_{z_i+j}$ is a nondecreasing parking function shifted by a
fixed value and $w_{z_i+1}\leq z_i+1$.
So $w$ is a parking function ($w_j\leq j$ for all $j$), and so is $uv$, since
it is a rearrangement of $w$.
This concludes the proof of the statement.

Note that $w$ is not in general nondecreasing but it satisfies nonetheless
$w_i\leq i$ for all $i$.
\qed
 
The lemma shows that whether $uv$ is a parking function or not depends only on
the biprofile of $(u,v)$, so that

\begin{corollary}
If $u$ and $v$ are nondecreasing words such that $uv$ is a parking function,
then any pair $(u',v')$ of nondecreasing words with the same biprofile as
$(u,v)$ is also such that $u'v'$ is a parking function.
\end{corollary}

\medskip
{\footnotesize
For example, consider the biprofile
\begin{equation}
\label{bipro-ex}
\binom{2\,5}{2\, 2} \ \text{and}\ \binom{1}{3}.
\end{equation}
There are four different choices for $u$: $2255$, $2256$, $2355$, and
$2356$ and five choices for $v$: $111$, $112$, $113$, $122$, and $123$.
One can check that all $20$ crossed concatenations are parking functions:
write down any $w=v.u$ and observe that even if $w$ is not always weakly
increasing, $w_i\leq i$ for all $i$.
}
\subsection{Regrouping terms in $\Delta G_n$}

We can now regroup the terms $\P^u\otimes \P^v$ in the ``unparkized''
multiplicity-free expression of $\Delta G_n$ according to their biprofiles,
and write
\begin{equation}
\Delta G_n =
  \sum_{\binom{s}{c},\binom{t}{d}}
    \sum_{\gf{\pf(u)=\binom{s}{c}}{\pf(v)=\binom{t}{d}}}
          \P^u\otimes \P^v
\end{equation}
where the sum runs over the parking biprofiles.

Now, given a parking biprofile $\binom{s}{c}, \binom{t}{d}$ and identifying
$w$ with $\Park(w)$, each sum

\begin{equation}
\sum_{\gf{\pf(u)
 = \binom{s}{c}}{\pf(v)=\binom{t}{d}}} \P^{\Park(u)} \otimes \P^{\Park(v)}
\end{equation}
contributes exactly one term $G^{c} \otimes G^{d}$ to $\Delta G$.
Indeed, given the profile $\binom{s}{c}$, the list of nondecreasing words
having that profile gives as parkized exactly all parking functions of size
$c_1$ contatenated with all parking functions of size $c_2$, \emph{etc.},
so that we get a term $G^c$.

\medskip
{\footnotesize
Continuing the example from Equation~\eqref{bipro-ex}, the term corresponding
to its biprofile is $G^{22}\otimes G^3$.
}

\medskip
Finally,
\begin{theorem}
$\Delta G_n$ is the sum of all $G^I\otimes G^J$, where $I$, $J$ run over
the bottom elements of all pairs $(b_1,b_2)$ of parking biprofiles of size
$n$.
\end{theorem}
\qed

Note in particular that if one swaps the profiles, one still has a parking
biprofile, which reflects the fact that $\Delta$ is cocommutative.

We shall also represent a profile as the minimal lexicographic nondecreasing
word associated with it. 

\medskip
{\footnotesize
For example, the biprofile 
$\binom{2\, 6\, 9}{3\,2\,2}$ is now $2226699$.
With this notation, the parking biprofiles of size $3$ correspond to the
following pairs of words:
\begin{equation}
\begin{split}
& (111,\emptyset), (11,1), (11,2), (11,3), (22,1), (13,1), (13,2), \\
& (1,11), (2,11), (3,11), (1,22), (1,13), (2,13), (\emptyset,111).
\end{split}
\end{equation}
so that sending a word to its packed evaluation,
\begin{equation}
\label{delta-g3}
\begin{split}
\Delta G_3 =
  G_3 \otimes 1 +
  (4\, G_2 + 2\, G^{11}) \otimes G_1 + 
  G_1 \otimes (4 \,G_2 + 2\, G^{11}) +
  1 \otimes G_3.
\end{split}
\end{equation}
}

\begin{note}
\label{deltaG-Catalan}
{\rm
Note that the number of terms in $\Delta G_3$ is $C_4=14$ (Catalan numbers),
and in general $\Delta G_n$ has $C_{n+1}$ terms.
A first easy but not very satisfactory goes as follows:
since $\Delta G_n$ is a sum of positive terms and that each $G_n$ is sent to
the usual $g_n$ when taking the commutative image from $\NCSF$ to $Sym$, each
term gives rise to one term of $\Delta g_n$. Since this coproduct is known to
have Catalan terms, so does $\Delta G_n$.

We will provide a complete combinatorial proof of this same result in the
Appendix through a bijection between parking biprofiles, pairs of
''compatible'' compositions and then Motzkin paths. It is also possible to
make a simple bijection between pairs of compositions and nondecreasing
parking functions but since this bijection does not provide any combinatorial
insight, we will only sketch it (see Note~\ref{comps-ndpf}).
}
\end{note}

\section{Combinatorial interpretations of $\Delta G_n$}

In the commutative case, it is known \cite{Einz} that 
\begin{equation}\label{eq:comcoprod}
\Delta g_n = \sum_{\pi\in\NC_{n+1}}g^{\alpha(\pi)}\otimes g^{\alpha(K(\pi))}
\end{equation}
where $\alpha(\pi)$ is the reduced type of $\pi$. We shall now see that this
expression can be extended to the noncommutative case, replacing the type by
the ordered type.

\subsection{From parking biprofiles to pairs of compositions}

We have seen that the coproduct of $g_n$ (or $G_n$, its pre-image in $\CQSym$)
can be expanded in the basis $g^I\otimes g^J$ and that the terms are
parametrized by parking biprofiles.

We shall now encode a profile $p$ by an integer composition $I$.

\begin{definition}
Let $p=\binom{s_1\dots s_k}{c_1\dots c_k}$ be a profile, and let 
$n\geq s_k+c_k$. Define $C:\ p\mapsto I$  as follows:
\begin{itemize}
\item If $s_1=1$ then $I=(1+c_1,I')$ where $I'$ is the composition associated
with the profile $\binom{s'_2\dots s'_k}{c_2\dots c_k}$ where $s'_i=s_i-c_1-1$.
\item If $s_1\not=1$ then $I=(1,I')$ where $I'$ is the composition associated
with the profile $\binom{s'_1\dots s'_k}{c_1\dots c_k}$ where $s'_i=s_i-1$.
\end{itemize}
Then, define $C_n(p)$ as the composition of $n$ obtained by adding $n-s_k-c_k$
ones at the end of $I$.
\end{definition}

\medskip
{\footnotesize
For example, with $n=12$,
\begin{equation}
\begin{split}
C\left(\binom{269}{221}\right)
 &= 1,C\left(\binom{158}{221}\right)
  = 1,3,C\left(\binom{25}{21}\right)
  = 1,3,1,C\left(\binom{14}{21}\right) \\
 &= 1,3,1,3,C\left(\binom{1}{1}\right)
  = 1,3,1,3,2,
\end{split}
\end{equation}
and finally $C_{12}(p)=(1,3,1,3,2,1,1)$.
Similarly,
\begin{equation}
C_{10}\left(\binom{16}{31}\right) = (4,1,2,1,1,1).
\end{equation}
}
\medskip
\begin{note}{\rm
Thanks to Lemma~\ref{carac-profils}, we know that a profile satisfies
$s_{i+1}>s_i+c_i$. Thus, at each step of the previous algorithm,
the $s'_i$ are positive integers, so that one indeed gets an integer
composition in the end.

Moreover, before adding ones at the end of $I$, one easily checks that $I$ was
a composition of $s_k+c_k$ so that $I$ itself is a composition of $n$.
}
\end{note}

The map $C$ is easily inverted:

\begin{lemma}
Let $I=(i_1,\dots,i_k)$ be a composition of $n$.
Define a map $C'$ by
\begin{equation}
C'(i_1,\dots,i_k) = \binom{d_1\ \ \dots\ \ d_k}{i_1\!-\!1\, \dots i_k\!-\!1},
\end{equation}
removing the biletters when the bottom letter is $0$ and
where $d_j=1+i_1+\dots+i_{j-1}$.

Then $C'$ is the inverse map of $C$.
\end{lemma}

\medskip
{\footnotesize
For example, with $I=(4,1,2,1,1,1)$, one gets $D=(1,5,6,8,9,10)$, so that
$C'(I)=\binom{16}{31}$ and $n=4+1+2+1+1+1=10$.

It will be useful to represent $I$ as a sequence of dots separated by bars,
such as
\begin{equation}
(4,1,2,1,1,1) \Longleftrightarrow ....|.|..|.|.|.
\end{equation}
On this representation, one easily reads $C'(I)$, and also the
lexicographically minimal  word $u$ with profile $C'(I)$: write an integer
equal to the position of the beginning of the block on each dot that is not
immediately followed by a bar.
On the example, we get
\begin{equation}
....|.|..|.|.|.  \Longrightarrow 111.|.|6.|.|.|.
\end{equation}
which indeed encodes $\binom{16}{31}$ and also its minimal word $1116$.
}

\medskip
\Proof By induction on the number of biletters of $p$.
Let $I=C_n(p)$.
If the first part of $I$ is not $1$, then we had $s_1=1$, its number of
occurrences $c_1$ being exactly $i_1-1$ by definition, so $C'$ records the
correct biletter at the beginning of its image.
The inductive definitions of the $s'$ and the $d$ are shifted in the same way
from $I$ to $I'$, which ensures that $C'(I')$ is the remaining part of $p$ by
induction.

If $I$ begins with a $1$, then $s_1$ was not $1$ and $d_1=1$ appears through
$C'$ with a $0$ at its bottom, so the biletter does not appear in $C'(I)$.
As before, the $s'$ and the $d$ change in the same way from $I$ to $I'$, so
$C'(I')$ will translate as $p$ by induction.
\qed

Now, given a parking biprofile $p=\binom{s}{c}$, $q=\binom{t}{d}$, we map it
to a pair of compositions by computing $C_n(p)$ and $C_n(q)$ with
$n=1+c_1+\dots+c_k+d_1+\dots+d_\ell$. Let us also denote this map by $C$.
Note that condition~\eqref{eq-parking-bipro} ensures that $n$ is greater than
both $s_k+c_k$ and $t_\ell+c_\ell$ so the map is well-defined and we get two
compositions of $n$.

\medskip
{\footnotesize
For example, 
\begin{equation}
C\left(\binom{269}{221}, \binom{16}{31} \right) =
\left( (1,3,1,3,2), (4,1,2,1,1,1)\right).
\end{equation}
}

\begin{definition}
\label{def-compat}
A pair of compositions is \emph{compatible} if it is in the image of $C$, that
is, the image of a parking biprofile.
\end{definition}

\begin{note}{\rm
Both $I$ and $J$ are compositions of the same integer $n$. Moreover, the
number of parts of $I$ is $n-(c_1+\dots+c_k)$ whereas the number of $J$ is
accordingly $n-(d_1+\dots+d_\ell)$, so that their total number of parts is
$n+1$.

Not all pairs satisfying this condition are compatible, but we shall see that
$I$ and its mirror conjugate $\bar I^\sim$ always are.
}
\end{note}

\begin{lemma}
A pair of compositions $(I,J)$ of the same integer $n$ is compatible iff their
total number of parts is $n+1$ and if the word $z$ obtained by sorting the
concatenation of the descent sets of $I$ and $J$ satisfies $z_\ell\geq \ell$
for all its values.
\end{lemma}

\medskip
{\footnotesize
For example, given the pair $I=(1,3,1,3,2)$ and $J=(4,1,2,1,1,1)$, the
concatenation of their descents is $[1,4,4,5,5,7,8,8,9]$ and it satisfies the
conditions of the statement.

As a counterexample, consider the pair $I=(1,2,1,1,2)$ and $J=(2,1,4)$. The
sorted concatenation of their descents is $[1,2,3,3,4,5]$ and $z_4=3<4$ so
that it does not satisfy the conditions of the statement and indeed,
$uv=261444$ has only two letters smaller than or equal to $3$.
}

\medskip
\Proof
We shall analyse the way in which the action of $C'$ on $uv$ depends on $z$.
To this aim, we shall represent a pair of compositions by  two sequences of
dots separated by bars. 

If $z_\ell\geq\ell$, there are at most $\ell-1$ bars among the $(\ell-1)$
first dots in the encodings of both $I$ and $J$. So there are at least
$2\ell-2-(\ell-1)=\ell-1$ values smaller than $\ell-1$ in $uv$, and the
parking constraint is fulfilled for $\ell-1$.

So, if $z_\ell\geq\ell$ for all $\ell$, then $C'(I,J)$ is a parking function.

Conversely, if some $z_\ell<\ell$, consider the smallest one $z_k$. 
Then,  $z_{k-1}$ was at least $k-1$ but since the word $z$ is weakly
increasing, we must have $z_{k-1}=z_k=k-1$.
In other words, both compositions $I$ and $J$ have a bar after $k-1$ dots and
there are also $k-2$ bars in total to the left of both these bars.
So among the $k-1$ dots on both lines of $I$ and $J$, exactly $(2k-2)-k=k-2$
do not have a dot immediately after them. Moreover, all the dots after the
$k$-th dot cannot be decodes as a value strictly smaller than $k$ since both
$I$ and $J$ have blocks beginning at position $k-1$. So there are exactly
$k-2$ values smaller than $k-1$ in $uv$ and so $uv$ is not a parking
function.
\qed

At this point, we have mapped bijectively the parking biprofiles to particular
pairs of compositions, and provided a characterization of those.
We can now use these results to provide an alternative description of the
coproduct of $G_n$.

\begin{lemma}
Through the bijection $C$, the map sending a parking biprofile to
$G^c\otimes G^d$ translates as the map sending a pair of compositions $(I,J)$
to $G^{i_1-1,\dots,i_r-1}\otimes G^{j_1-1,\dots,j_r-1}$
and removing the zeroes.
\end{lemma}

\Proof
Immediate by definition of $C$.
\qed

\medskip
{\footnotesize
Here follows the whole list of compatible pairs of compositions of size $4$:
\begin{equation}
\begin{split}
& (4,1111), (31,211), (31,121), (31,112), (22,211), (22,121), (211,31), \\
& (211,22), (211,13), (13,211), (121,31), (121,22), (112,31), (1111,4).
\end{split}
\end{equation}
and one can then check the expression of $G_3$ of~\eqref{delta-g3} by sending
each composition $I$ to $G^{i_1-1,\dots,i_r-1}$.
}

\begin{note}
\label{comps-ndpf}
	{\rm
We shall provide in the Appendix (Section~\ref{sec-appendix}) a meaningful
bijection proving that pairs of compositions are enumerated by Catalan numbers
but we can provide a very simple but dumb one that also proves that: given a
pair $(I,J)$ of $n$ of respective descent sets $(d_1,\dots,d_k)$ and
$(d'_1,\dots,d'_\ell)$, sort the word $w$ given by the concatenation of
the $2*d_i-1$ with the $2*d'_i$ and the value $2n-1$.
Now compute $w'$ where $w'_{n+1-i}=n+i-w_i$.

This is a bijection from the pairs of compositions to their image set since
both set are easily revertible. And it is an exercice to check that $w'$ is a
nondecreasing parking function and conversely that any parking function gives
rise to a valid pair of compositions.
	}
\end{note}

\medskip
{\footnotesize
Given the pair $(13132,412111)$, one gets the descents sets $(1,4,5,8)$ and
$(4,5,7,8,9)$ hence the word
\begin{equation}
w = 1,7,8,9,10,14,15,16,18,19
\text{\ \ \ and\ \ \ }
w' = 1, 1, 2, 2, 2, 5, 5, 5, 5, 10
\end{equation}
that is indeed a nondecreasing parking function.
}

\subsection{From pairs of compositions to noncrossing partitions}

\subsubsection{Noncrossing partitions and permutations}

Recall that a noncrossing partition $\pi$ can be interpreted as a permutation
$w_\pi$ whose cycles are the blocks of $\pi$ read in increasing order. The
\emph{right Kreweras complement}~\cite{Kr}  $\pi'=K(\pi)$ can then be defined
as the noncrossing partition such that $w_{\pi'}=w_\pi^{-1}\gamma_n$, where
$\gamma_n=(123\dots n)$ is the canonical long cycle.
The permutations $w_\pi$ are called \emph{noncrossing permutations}.

\medskip
{\footnotesize
For example, given the noncrossing partition $\pi$,
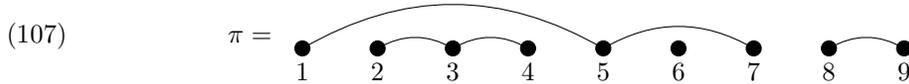
\begin{equation}
\pi= \raisebox{-2em}{\Matching{9}{1/5, 2/3, 3/4, 5/7,8/9}}
\end{equation}
we get
\begin{equation}
w_\pi = [(1,5,7)(2,3,4)(8,9)] = [5,3,4,2,7,6,1,9,8]
\end{equation}
so that
\begin{equation}
\begin{split}
w_{\pi'} &= [7,4,2,3,1,6,5,9,8].[2,3,4,5,6,7,8,9,1]  \\
         &= [4,2,3,1,6,5,9,8,7] = [(1,4)(5,6)(7,9)]
\end{split}
\end{equation}
and $\pi'$ is
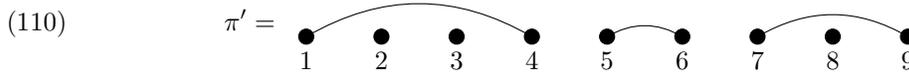
\begin{equation}
\pi'= \raisebox{-2em}{\Matching{9}{1/4, 5/6, 7/9}}
\end{equation}
}

\medskip
The \emph{canonical ordering} of a permutation is the list of its cycles
in increasing order of their minimal elements.

\subsubsection{Planar binary trees and the Kreweras complement}

There are many bijections between noncrossing partitions and binary trees.
But actually, on a binary tree $t$, one can directly read \emph{two}
noncrossing partitions $\pi'$, $\pi''$. 

Let $\treetoncpncp$ be the map sending a tree $T$ to a pair $(\pi',\pi'')$
as follows : label the nodes of $T$ in infix order, so as to obtain a binary
search tree:

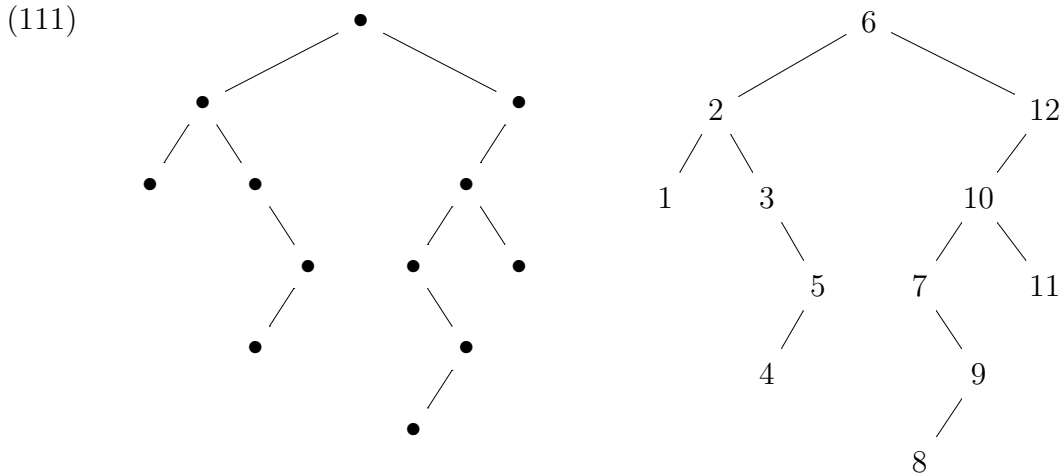
\begin{equation}
\label{arb-12}
\setlength\unitlength{1.7mm}
\xymatrix@R=0.5cm@C=.9mm{
 & & & &
   *{\GrTeXBox{\bullet}}\arx1[llld]\arx1[rrrd] & \\
 & *{\GrTeXBox{\bullet}}\arx1[ld]\arx1[rd] & & & & & & 
   *{\GrTeXBox{\bullet}}\arx1[ld] & \\
   *{\GrTeXBox{\bullet}} & &
   *{\GrTeXBox{\bullet}}\arx1[rd] & & & &
   *{\GrTeXBox{\bullet}}\arx1[ld]\arx1[rd] \\
 & & &
   *{\GrTeXBox{\bullet}}\arx1[ld] & & 
   *{\GrTeXBox{\bullet}}\arx1[rd] & &
   *{\GrTeXBox{\bullet}} \\
 & &
   *{\GrTeXBox{\bullet}} & & & & 
   *{\GrTeXBox{\bullet}}\arx1[ld] \\
 & & & & & 
   *{\GrTeXBox{\bullet}} \\
}
\hskip1cm
\xymatrix@R=0.5cm@C=.7mm{
 & & & &
   *{\GrTeXBox{6}}\arx1[llld]\arx1[rrrd] & \\
 & *{\GrTeXBox{2}}\arx1[ld]\arx1[rd] & & & & & & 
   *{\GrTeXBox{12}}\arx1[ld] & \\
   *{\GrTeXBox{1}} & &
   *{\GrTeXBox{3}}\arx1[rd] & & & &
   *{\GrTeXBox{10}}\arx1[ld]\arx1[rd] \\
 & & &
   *{\GrTeXBox{5}}\arx1[ld] & & 
   *{\GrTeXBox{7}}\arx1[rd] & &
   *{\GrTeXBox{11}} \\
 & &
   *{\GrTeXBox{4}} & & & & 
   *{\GrTeXBox{9}}\arx1[ld] \\
 & & & & & 
   *{\GrTeXBox{8}} \\
}
\end{equation}

Then the blocks of $\pi'$ are the sets of labels of the left branches of $T$:
{\footnotesize
\begin{equation}
\pi'= \raisebox{-2em}{\Matching{12}{1/2,2/6,4/5,7/10,10/12,8/9}}
\end{equation}
}
and the blocks of $\pi''$ are those of its right branches:
{\footnotesize
\begin{equation}
\pi''= \raisebox{-2em}{\Matching{12}{2/3,3/5,6/12,7/9,10/11}}.
\end{equation}
}
Both $\pi'$ and $\pi''$ are obviously noncrossing partitions.
Moreover, if one traverses the tree in infix order and records the labels of
each branch the first time it is encountered (that is, by its smallest value),
both partitions $\pi'$ and $\pi''$ come up with their canonical ordering.
It is also easy to see that $\treetoncpncp$ is also bijective since one can
easily rebuild $T$ from either $\pi'$ or $\pi''$. This means that one of
the elements should be fully recoverable from the other, or, in other words,
that they have a direct link with one another, and indeed:

\medskip
{\footnotesize
Interpreting $\pi'$ and $\pi''$ as permutations,
\begin{equation}
\pi'
 = [(126), (3), (4, 5), (7, 10, 12), (8, 9), (11)]
 = [2, 6, 3, 5, 4, 1, 10, 9, 8, 12, 11, 7]
\end{equation}
and
\begin{equation}
\pi''
 = [(1), (2, 3, 5), (4), (6, 12), (7, 9), (8), (10, 11)]
 = [1, 3, 5, 4, 2, 12, 9, 8, 7, 11, 10, 6]
\end{equation}
and one can check that
\begin{equation}
\pi'\pi''=[2, 3, 4, 5, 6, 7, 8, 9, 10, 11, 12, 1]
\end{equation}
so that $\pi''$ is the right Kreweras complement of $\pi'$.
}

\medskip
It is easy to see that this is true in general:

\begin{lemma}
\label{lem-gd}
Let $T$ be a binary tree and let $\treetoncpncp(T)=(\pi',\pi'')$.
Then $\pi''$ is the right Kreweras complement of $\pi'$.
\end{lemma}

\Proof
The property holds for trees with at most $2$ nodes and also for trees with no
right or left branches since in these cases, either $\pi'$ or $\pi''$ is
the cycle $\gamma_n$ and the other is the identity.

Assume by induction that the property holds for trees with at most $n-1$
nodes. Let $T$ be a tree whose left subtree $T_L$ has $k-1$ nodes and whose
right subtree $T_R$ has $n-k$ nodes. Its root has therefore label $k$.
By induction hypothesis (or the special case mentioned above), the product of
the cycles of the tree $T'_L$ with root $k$, left subtree $T_L$ and an empty
right subtree is the cycle $\sigma'_1=(1\cdots k-1k)$.
Similarly, the product of the cycles of the tree $T'_R$ with root $k$, empty
left subtree and right subtree $T_R$ is $\sigma'_2=(k\ k+1\cdots n)$. The
complete product is therefore $\sigma'_1\sigma'_2=(12\cdots n)$.
\qed

\begin{note}{\rm
\label{note-infix}
Lemma~\ref{lem-gd} amounts to saying that, inside a binary search tree, one
gets from the position of $i$ to the position of $i+1$ modulo its number of
nodes by
\begin{itemize}
\item moving one step down its right branch (and cycling if $i$ is at the
bottom of it),
\item moving then one step up its left branch (and cycling if it is at
the top of it).
\end{itemize}
This property is easily checked, and provides an alternative proof of the
lemma.
}
\end{note}

\subsubsection{From trees and permutations to pairs of compositions}

As we have seen, reading the left branches and the right branches in infix
order, the blocks of the partitions come naturally ordered with respect to
their minima in increasing order, so that the compositions recording their
lengths are the ordered types of $\pi'$ and $\pi''$.

It turns out that the tree $T$, and therefore $\pi'$ and $\pi''$ can be
unambiguously reconstructed from this pair of compositions.

\begin{theorem}
\label{bicomp-inj}
Let $t(\pi)$ denote the ordered type of a noncrossing partition. The map
\begin{equation}
\tau: \pi\mapsto (t(\pi),t(K(\pi))
\end{equation}
is injective.
\end{theorem}

The map $\tau$ goes from a noncrossing partition to a pair of compositions.
Since one can easily go from a noncrossing partition to a tree, we shall also
write $\tau$ as the map sending a tree to a pair of compositions and it is
that map, sending a tree $T$ to the lengths of the ordered types of
$\treetoncpncp(T)$, that we will prove injective.

Let us consider the following backwards algorithm:

\begin{algorithm}
\label{algo-rebuild}
Input:  a pair of compositions $I=(i_1,\ldots,i_r)$ and $J=(j_1,\ldots,j_s)$
obtained as the ordered lengths of the respective left and right branches of a
tree.

We shall build a tree one branch at each step. When gluing a branch on a
node, mark this node.

Create a left branch of $i_1$ nodes.
Then glue a right branch of $j_1$ nodes at the first unmarked node in infix
order (in that case, it is the leftmost node since no one was marked yet).

Then move to the first unmarked node in infix order, which can be either a
left or a right child (if it is the root, consider it as a left child), and
create a new branch in the opposite direction (\emph{e.g.}, right if it is a
left child) of the corresponding size, using the next unused part of $I$ or of
$J$ depending on the direction.
Iterate until there are no unmarked nodes left.
\end{algorithm}

\medskip
{\footnotesize
An example of this algorithm with $I=312321$ and $J=1312212$ is given on
Fig.~\ref{fig-algo}. Note that when a part is $1$, we just mark the
leftmost node and add no new node.
}

\def\cpl{\circle*{1.2}}
\def\cvi{\circle{1.2}}
\def\cu{1.2}
\def\cv{1.4}
\def\cor{.6}

\def\mmd{-1.2}
\def\ppd{1.2}
\def\mmt{-2}
\def\ppt{2}

\def\aau{
\begin{picture}(4,2)
 \put(2,2){\circle*{\cu}}
\end{picture}}

\def\aad{\begin{picture}(8,5)
 \put(6,6){\cvi}
 \put(5.6,5.6){\Line(\mmd,\mmd)}
 \put(4,4){\cvi}
 \put(3.6,3.6){\Line(\mmd,\mmd)}
 \put(2,2){\cvi}
\end{picture}}

\def\aat{\begin{picture}(8,5)
 \put(6,6){\cvi}
 \put(5.6,5.6){\Line(\mmd,\mmd)}
 \put(4,4){\cvi}
 \put(3.6,3.6){\Line(\mmd,\mmd)}
 \put(2,2){\cpl}
\end{picture}}

\def\aaq{\begin{picture}(8,5)
 \put(6,6){\cvi}
 \put(5.6,5.6){\Line(\mmd,\mmd)}
 \put(4,4){\cpl}
 \put(3.6,3.6){\Line(\mmd,\mmd)}
 \put(4.4,3.6){\Line(\ppd,\mmd)}
 \put(2,2){\cpl}
 \put(6,2){\cvi}
 \put(6.4,1.6){\Line(\ppd,\mmd)}
 \put(8,0){\cvi}
\end{picture}}

\def\aac{\begin{picture}(8,5)
 \put(6,6){\cvi}
 \put(5.6,5.6){\Line(\mmd,\mmd)}
 \put(4,4){\cpl}
 \put(3.6,3.6){\Line(\mmd,\mmd)}
 \put(4.4,3.6){\Line(\ppd,\mmd)}
 \put(2,2){\cpl}
 \put(6,2){\cpl}
 \put(6.4,1.6){\Line(\ppd,\mmd)}
 \put(8,0){\cvi}
\end{picture}}

\def\aas{\begin{picture}(8,5)
 \put(6,6){\cvi}
 \put(5.6,5.6){\Line(\mmd,\mmd)}
 \put(4,4){\cpl}
 \put(3.6,3.6){\Line(\mmd,\mmd)}
 \put(4.4,3.6){\Line(\ppd,\mmd)}
 \put(2,2){\cpl}
 \put(6,2){\cpl}
 \put(6.4,1.6){\Line(\ppd,\mmd)}
 \put(8,0){\cpl}
 \put(7.6,-0.4){\Line(\mmd,\mmd)}
 \put(6,-2){\cvi}
\end{picture}}

\def\aae{\begin{picture}(8,5)
 \put(6,6){\cvi}
 \put(5.6,5.6){\Line(\mmd,\mmd)}
 \put(4,4){\cpl}
 \put(3.6,3.6){\Line(\mmd,\mmd)}
 \put(4.4,3.6){\Line(\ppd,\mmd)}
 \put(2,2){\cpl}
 \put(6,2){\cpl}
 \put(6.4,1.6){\Line(\ppd,\mmd)}
 \put(8,0){\cpl}
 \put(7.6,-.4){\Line(\mmd,\mmd)}
 \put(6,-2){\cpl}
\end{picture}}

\def\aah{\begin{picture}(16,5)
 \put(4,4){\cpl}
 \put(3.6,3.6){\Line(\mmd,\mmd)}
 \put(4.4,3.6){\Line(\ppd,\mmd)}
 \put(2,2){\cpl}
 \put(6,2){\cpl}
 \put(6.4,1.6){\Line(\ppd,\mmd)}
 \put(8,0){\cpl}
 \put(7.6,-.4){\Line(\mmd,\mmd)}
 \put(6,-2){\cpl}
 \put(10,8){\cpl}
 \put(10,8){\Line(-6,-4)}
 \put(10,8){\Line( 5.5,-3.6)}
 \put(16,4){\cvi}
\end{picture}}

\def\aan{\begin{picture}(16,5)
 \put(10,8){\cpl}
 \put(10,8){\Line(-6,-4)}
 \put(10,8){\Line( 5.5,-3.6)}
 \put(4,4){\cpl}
 \put(3.6,3.6){\Line(\mmd,\mmd)}
 \put(4.4,3.6){\Line(\ppd,\mmd)}
 \put(2,2){\cpl}
 \put(6,2){\cpl}
 \put(6.4,1.6){\Line(\ppd,\mmd)}
 \put(8,0){\cpl}
 \put(7.6,-.4){\Line(\mmd,\mmd)}
 \put(6,-2){\cpl}
 \put(16,4){\cpl}
 \put(15.6,3.6){\Line(\mmd,\mmd)}
 \put(14,2){\cvi}
 \put(13.6,1.6){\Line(\mmd,\mmd)}
 \put(12,0){\cvi}
\end{picture}}

\def\aadi{\begin{picture}(16,5)
 \put(10,8){\cpl}
 \put(10,8){\Line(-6,-4)}
 \put(10,8){\Line( 5.5,-3.6)}
 \put(4,4){\cpl}
 \put(3.6,3.6){\Line(\mmd,\mmd)}
 \put(4.4,3.6){\Line(\ppd,\mmd)}
 \put(2,2){\cpl}
 \put(6,2){\cpl}
 \put(6.4,1.6){\Line(\ppd,\mmd)}
 \put(8,0){\cpl}
 \put(7.6,-.4){\Line(\mmd,\mmd)}
 \put(6,-2){\cpl}
 \put(16,4){\cpl}
 \put(15.6,3.6){\Line(\mmd,\mmd)}
 \put(14,2){\cvi}
 \put(13.6,1.6){\Line(\mmd,\mmd)}
 \put(12,0){\cpl}
 \put(12.4,-.4){\Line(\ppd,\mmd)}
 \put(14,-2){\cvi}
\end{picture}}

\def\aaon{\begin{picture}(16,5)
 \put(10,8){\cpl}
 \put(10,8){\Line(-6,-4)}
 \put(10,8){\Line( 5.5,-3.6)}
 \put(4,4){\cpl}
 \put(3.6,3.6){\Line(\mmd,\mmd)}
 \put(4.4,3.6){\Line(\ppd,\mmd)}
 \put(2,2){\cpl}
 \put(6,2){\cpl}
 \put(6.4,1.6){\Line(\ppd,\mmd)}
 \put(8,0){\cpl}
 \put(7.6,-.4){\Line(\mmd,\mmd)}
 \put(6,-2){\cpl}
 \put(16,4){\cpl}
 \put(15.6,3.6){\Line(\mmd,\mmd)}
 \put(14,2){\cvi}
 \put(13.6,1.6){\Line(\mmd,\mmd)}
 \put(12,0){\cpl}
 \put(12.4,-.4){\Line(\ppd,\mmd)}
 \put(14,-2){\cpl}
 \put(13.6,-2.4){\Line(\mmd,\mmd)}
 \put(12,-4){\cvi}
\end{picture}}

\def\aado{\begin{picture}(16,5)
 \put(10,8){\cpl}
 \put(10,8){\Line(-6,-4)}
 \put(10,8){\Line( 5.5,-3.6)}
 \put(4,4){\cpl}
 \put(3.6,3.6){\Line(\mmd,\mmd)}
 \put(4.4,3.6){\Line(\ppd,\mmd)}
 \put(2,2){\cpl}
 \put(6,2){\cpl}
 \put(6.4,1.6){\Line(\ppd,\mmd)}
 \put(8,0){\cpl}
 \put(7.6,-.4){\Line(\mmd,\mmd)}
 \put(6,-2){\cpl}
 \put(16,4){\cpl}
 \put(15.6,3.6){\Line(\mmd,\mmd)}
 \put(14,2){\cvi}
 \put(13.6,1.6){\Line(\mmd,\mmd)}
 \put(12,0){\cpl}
 \put(12.4,-.4){\Line(\ppd,\mmd)}
 \put(14,-2){\cpl}
 \put(13.6,-2.4){\Line(\mmd,\mmd)}
 \put(12,-4){\cpl}
\end{picture}}

\def\aatr{\begin{picture}(16,5)
 \put(10,8){\cpl}
 \put(10,8){\Line(-6,-4)}
 \put(10,8){\Line( 5.5,-3.6)}
 \put(4,4){\cpl}
 \put(3.6,3.6){\Line(\mmd,\mmd)}
 \put(4.4,3.6){\Line(\ppd,\mmd)}
 \put(2,2){\cpl}
 \put(6,2){\cpl}
 \put(6.4,1.6){\Line(\ppd,\mmd)}
 \put(8,0){\cpl}
 \put(7.6,-.4){\Line(\mmd,\mmd)}
 \put(6,-2){\cpl}
 \put(16,4){\cpl}
 \put(15.6,3.6){\Line(\mmd,\mmd)}
 \put(14,2){\cpl}
 \put(13.6,1.6){\Line(\mmd,\mmd)}
 \put(12,0){\cpl}
 \put(12.4,-.4){\Line(\ppd,\mmd)}
 \put(14,-2){\cpl}
 \put(13.6,-2.4){\Line(\mmd,\mmd)}
 \put(12,-4){\cpl}
 \put(14.4,1.6){\Line(\ppd,\mmd)}
 \put(16,0 ){\cvi}
\end{picture}}

\def\aaqu{\begin{picture}(8,8)
 \put(10,8){\cpl}
 \put(10,8){\Line(-6,-4)}
 \put(10,8){\Line( 5.5,-3.6)}
 \put(4,4){\cpl}
 \put(3.6,3.6){\Line(\mmd,\mmd)}
 \put(4.4,3.6){\Line(\ppd,\mmd)}
 \put(2,2){\cpl}
 \put(6,2){\cpl}
 \put(6.4,1.6){\Line(\ppd,\mmd)}
 \put(8,0){\cpl}
 \put(7.6,-.4){\Line(\mmd,\mmd)}
 \put(6,-2){\cpl}
 \put(16,4){\cpl}
 \put(15.6,3.6){\Line(\mmd,\mmd)}
 \put(14,2){\cpl}
 \put(13.6,1.6){\Line(\mmd,\mmd)}
 \put(12,0){\cpl}
 \put(12.4,-.4){\Line(\ppd,\mmd)}
 \put(14,-2){\cpl}
 \put(13.6,-2.4){\Line(\mmd,\mmd)}
 \put(12,-4){\cpl}
 \put(14.4,1.6){\Line(\ppd,\mmd)}
 \put(16,0 ){\cpl}
\end{picture}}

\begin{figure}[ht]
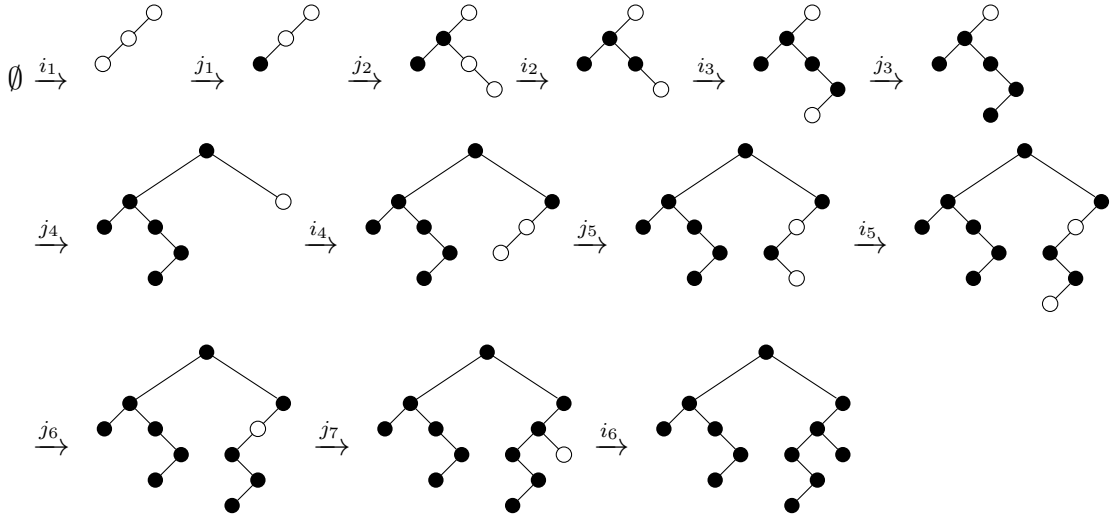

\begin{equation*}
\begin{split}
\emptyset
     &\xrightarrow{i_1} \aad
      \xrightarrow{j_1} \aat
      \xrightarrow{j_2} \aaq \
      \xrightarrow{i_2} \aac \ \
      \xrightarrow{i_3} \aas \ \
      \xrightarrow{j_3} \aae \ \\[1cm]
    & \xrightarrow{j_4} \aah \ 
      \xrightarrow{i_4} \aan \
      \xrightarrow{j_5} \aadi\ \
      \xrightarrow{i_5} \aaon\ \\[1cm]
    & \xrightarrow{j_6} \aado\ \
      \xrightarrow{j_7} \aatr\ \
      \xrightarrow{i_6} \aaqu\ \\[.5cm]
\end{split}
\end{equation*}
\caption{\label{fig-algo}Algorithm~\ref{algo-rebuild} applied to the image of
the tree in Equation~\eqref{arb-12}. Unmarked nodes are white.}
\end{figure}

\begin{proposition}
\label{prop-algoinv}
If $I$ and $J$ are the ordered lengths of the left and right branches of a
tree $T$, then Algorithm~\ref{algo-rebuild} rebuilds $T$ from $I$ and $J$.
\end{proposition}

\Proof
Let $T$ be a tree and let $\tau(T)=(I,J)$.

Apply Algorithm~\ref{algo-rebuild} to $I$ and $J$.
We shall prove by induction that after step $k$ the partial tree is
exactly the tree $T'$ consisting of the first $k$ left and/or right branches
of $T$.

This is true at steps $k=0$ and $k=1$. Assume that this is true until step $k$
and add a new (left or right) branch to $T'$ as described in
Algorithm~\ref{algo-rebuild}.
Let  $T''$ be the resulting tree. Note that if $T'$ is not equal to $T$, there
are necessarily unmarked nodes so that step $k+1$ is well-defined.

By construction, this branch has been added to the leftmost unmarked node $x$
in infix order. By definition of this order, all marked nodes strictly before
$x$ will all be read in the same order as in $T'$, all read before $x$ and its
added branch. So the first $k+1$ branches associated with $T'$ are the first
$k$ branches of $T$ followed by the added new branch.
\qed

\Proof[of the theorem]
Thanks to Lemma~\ref{prop-algoinv}, the map $\tau$  from binary trees to
pairs of compositions induces a bijection with  its image set.
\qed 

We finally need to characterize the image set of $\tau$.

\begin{lemma}
Let $T$ be a tree and $\tau(T)=(I,J)$.

Then $I$ and $J$ are compatible in the sense of Definition~\ref{def-compat}.
\end{lemma}

\Proof
First of all, it is well-known that if $\pi'$ and $\pi''$ satisfy
$\pi'\pi''=\gamma_n$ and $\ell(\pi')+\ell(\pi'')=n$, then their total number
of cycles is $n+1$.
So $I$ and $J$, being the images of two such permutations have a total of $n$
descents.

Now, let us consider two compositions $I$ and $J$ whose total number of
descents is $n$. Sort these descents, and let $d=d_1\dots d_n$ be the
corresponding word. Either they satisfy the criterion or there is a smallest
value $k$ such that $d_k<k$.
Since $d_{k-1}\geq k-1$, both $d_{k-1}=d_k=k-1$ and $I$ and $J$ both have a
descent in $k-1$.

Let $I'=(i_1,\dots,i_{k_1})$ and $J'=(j_1,\dots,j_{k_2})$ be the prefixes of
$I$ and $J$ such that $i_1+\dots+i_{k_1}=j_1+\dots+j_{k_2}=k-1$. Then $I'$ and
$J'$ are both compositions of $k-1$ whose total number of descents is $k-2$,
since we do not take into account their last descent. Moreover, their own
descents are the first $d_k$s so $I'$ and $J'$ are compatible. By induction,
they correspond therefore to a tree $T'$.

Let us no apply Algorithm~\ref{algo-rebuild} to $I$ and $J$. Since the
algorithm proceeds step by step, if it used part $i_{k_1+1}$ or part
$j_{k_2+1}$ before going through $I'$ and $J'$ fully,  the algorithm would
have failed on the pair $(I',J')$, which is not the case. 
So it ends with $T'$ at this exact step, there is no unmarked node left and
the algorithm stops.

So, by induction, the algorithm fails if $I$ and $J$ are not compatible. Since
we know that the algorithm succeeds with the images of the binary trees, it
means that the image set of $\tau$ is included in the set of compatible pairs
of compositions. But both sets are Catalan sets (see
Notes~\ref{deltaG-Catalan} and~\ref{comps-ndpf}) so they coincide.
\qed

\subsubsection{Conclusion of the proof of Theorem \ref{th:coprod}}

We have succesfully expressed $\Delta G_n$ as a sum over parking biprofiles
$G^c\otimes G^d$, then mapped parking biprofiles to compatible pairs of
compositions, and proved that such pairs record the ordered lengths of left
and right branches of binary trees. Such pairs of compositions in turn
coincide with the reduced types of $(\pi,K(\pi))$ for $\pi\in\NC_{n+1}$.
This concludes the proof of Theorem~\ref{th:coprod}.
\qed

\begin{note}{\rm
On the interpretation of $\Delta g_n$ in terms of noncrossing partitions and
their  Kreweras complement, it is not apparent that $\Delta$ is cocommutative, since this
operation is not an involution. It would then be interesting, given a
noncrossing partition $\pi$ of $[n+1]$ of reduced ordered type $I$ whose
Kreweras complement $\pi'$ has reduced ordered type $J$, to build a
noncrossing partition $\pi'$ of $[n+1]$ of reduced ordered type $J$ whose
Kreweras complement $\pi'$ has reduced ordered type $I$.

The known involutions on noncrossing partitions, iterations of Kreweras and
the same up to reversal defined by Simion and Ullman in~\cite{SU}, do not
have this property.
}
\end{note}

{\footnotesize
For example, with $I=(4,1,1,1,2,2,1)$ and $J=(3,1,3,1,2,2)$,
the map would exchange
\begin{equation}
p := [4, 2, 3, 5, 9, 6, 8, 7, 1, 12, 11, 10]
\text{\ \ and\ \ }
p' := [7, 2, 4, 5, 3, 6, 8, 1, 12, 11, 10, 9]
\end{equation}
of respective cycles
\begin{equation}
c_0  = [(1, 4, 5, 9), (2), (3), (6), (7, 8), (10, 12), (11)]
\text{\ and\ }
c'_0 = [(1, 7, 8), (2), (3, 4, 5), (6), (9, 12), (10, 11)]
\end{equation}
whose  Kreweras complements have as cycles
\begin{equation}
c_1  = [(1, 2, 3), (4), (5, 6, 8), (7), (9, 12), (10, 11)]
\text{\ and\ }
c'_1 = [(1, 2, 5, 6), (3), (4), (7), (8, 12), (9, 11), (10)].
\end{equation}
}

\section{Coproduct of $g_n$ in the commutative case}

The calculation of $\Delta g_n$ is easier in the commutative case.
We have seen in Lemma~\ref{lem-gd} that one can read a noncrossing
partition $\pi$ and its Kreweras complement on a binary tree. 
This information is enough to compute the coproduct of $g(X)$.

Indeed, $g(X+Y)$ satisfies the functional equation
\begin{align}
g(X+Y)&= \sigma_{g(X+Y)}(X+Y)=\sigma_{g(X+Y)}(X)\sigma_{g(X+Y)}(Y)\\
& = \sum_{_p\ge 0}h_p(X)g(X+Y)^p\sum_{q\ge 0}h_q(Y)g(X+Y)^q,
\end{align}
\emph{i.e.}, the right-hand side factorizes. This is not true anymore in the
noncommutative case, and we had to rely upon a different argument, based on
the possibility to reconstruct $\pi$ and $K(\pi)$ from their ordered types.

In the commutative case, we shall show that this equation coincides with that
of the generating series of binary trees by lengths of the left and right
branches. The argument is similar to the one used by Goulden and Jackson in
their proof of Macdonald's formula for the top connexion
coefficients~\cite{GJ1,GJ2}. This provides an alternative (and simpler) proof
of the result of~\cite{Einz}.

For a binary tree $t$, set
\begin{equation}
w(t;u,v) = \prod_{\ell\in L(t)} u_{e(\ell)} \prod_{r\in R(t)}v_{e(r)}\,,
\end{equation}
where $L(t)$ and $R(t)$ are respectively the sets of left and right branches
of $t$, and $e(b)$ denotes the number of edges in a branch $b$ (with the
convention $u_0=v_0=1$).

\medskip
{\footnotesize
For example, we have $w(t;u,v)=u_2^2u_1^2\cdot v_2v_1^3$ on the
following tree
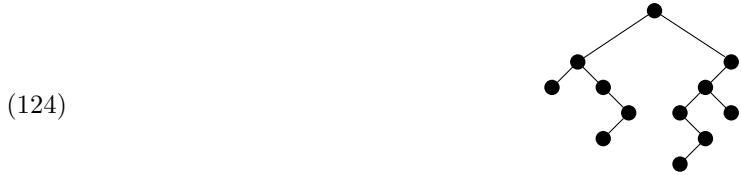
\begin{equation}
\aaqu
\end{equation}
\bigskip
\medskip
}

Let $W$ be the generating series
\begin{equation}
W(u,v)=\sum_{t\in\BT}w(t;u,v) = 1 + u_1+v_1 + u_2+3u_1v_1+v_2 + \cdots
\end{equation}
Denote by $t_L$ and  $t_R$ the left and right subtrees of $t$, and let $U,V$
be the generating series of the trees whose right (resp. left) subtree is
empty:
\begin{equation}
U = \sum_{t_R=\emptyset}w(t;u,v),\quad V = \sum_{t_L=\emptyset}w(t;u,v).
\end{equation}
These series satisfy the system
\begin{equation}
\label{eq:system}
\begin{cases}
V &= \displaystyle\sum_{n\ge 0}v_n U^n\\
\medskip
U &= \displaystyle\sum_{n\ge 0}u_nV^n,
\end{cases}
\end{equation}
and classifying trees by length of the left branch of the root, we can write
\begin{equation}
W = V +u_1V^2+u_2V^3+u_3V^4+\cdots = UV.
\end{equation}
Recall that the commutative symmetric Lagrange series solves the equation
\begin{equation}
t = \frac{u}{\sigma_u(X)}\Longleftrightarrow u=tg(t;X)
  = \sum_{n\ge 0}g_n(X)t^{n+1}.
\end{equation}
If in \eqref{eq:system} we set $u_n=g_n(X)$ and $v_n=g_n(Y)$, the system
becomes, multiplying the first equation by $U$ and the second one by $V$

\begin{equation}
\label{eq:system2}
\begin{cases}
UV &= \displaystyle\sum_{n\ge 0}g_n(Y) U^{n+1}\\
\medskip
UV &= \displaystyle\sum_{n\ge 0}g_n(X)V^{n+1},
\end{cases}
\end{equation}
whence
\begin{equation}
UV = Ug(U;Y)\Leftrightarrow U 
   = \frac{UV}{\sigma_{UV}(Y)}\Leftrightarrow V
   = \sigma_{UV}(Y)
\end{equation}
and similarly
\begin{equation}
UV = Vg(U;X)\Leftrightarrow V 
   = \frac{UV}{\sigma_{UV}(X)}\Leftrightarrow U
   = \sigma_{UV}(X)
\end{equation}
Therefore,

\begin{equation}
W = UV = \sigma_{UV}(X+Y)
\end{equation}
which is precisely the functional equation of $\Delta g$.

This proves Equation \eqref{eq:comcoprod}.

\section{Application to the reduced incidence Hopf algebra of noncrossing
partitions}

In the commutative case, the calculation of $\Delta g_n$ and $g_n(-X)$,
of which we have given new proofs, have important applications to the
combinatorics of noncrossing partitions.
Although the results of this Section are known, it seems appropriate to take
the opportunity of giving a streamlined account of the theory in the light of
the previous considerations.

The reduced incidence Hopf algebra $\Hh_\NC$ of the hereditary family of
lattices $\NC_n$ is the vector space spanned by isomorphism classes of
intervals of the $\NC_n$ for $n\ge 1$. The order is defined by $\pi\le\pi'$
if $\pi$ is finer than $\pi'$. The minimal element $\O_n$ is the partition
into singletons, and the maximal element $\I_n$ is the partition with one
block. As is well-known \cite{Spei1}, any such interval is isomorphic to a
Cartesian product of complete lattices $\NC_k$.
An interval $[\O_n, \pi]$ is isomorphic to $\prod_{B\in \pi}\NC_{|B|}$, and an
interval $[\pi,\I_n]$ is isomorphic to $[\O_n,K(\pi)]$. Finally, if
$\sigma=\{B_1,\ldots,B_r\}$, $[\pi,\sigma]$ is
isomorphic to $\prod_i[\pi\cap B_i,\I_{B_i}]$.

The product of $\Hh_\NC$ is the Cartesian product. Thus, $\Hh_\NC$ is the
polynomial algebra on the variables $y_n=[\NC_{n+1}]$, and the coproduct is
defined as
\begin{equation}\label{eq:dy}
	\delta y_n = \sum_{\pi\in\NC_{n+1}}[\O_{n+1},\pi]\otimes [\pi,\I_{n+1}].
\end{equation}
One of the main results of \cite{Einz} shows that $y_n\mapsto g_n(-X)$ is an
isomorphism of Hopf algebras from $\Hh_{\NC}$ to $Sym$. This is precisely what
we have just proved (using $g_n(X)$ instead) by a different method.

Another result of \cite{Einz}, which has been reproved by a different method
in \cite{EH} is equivalent to the computation of the antipode $g_n(-X)$ in
$Sym$. It is implied by our calculation of $g(-A)$ (in \cite{Einz} the
coefficients are interpreted as counting polygon dissections, but our formula
is cancellation-free as well, and produces the same coefficients).

It should be noted that these calculations imply a great deal of classical
results about noncrossing partitions. In particular, the multiplicative
functions on noncrossing partitions are the characters of $\Hh_\NC$. Such a
function $\phi$ is completely determined by its values $a_n=\phi(y_n)$ on the
generators.

Using the above isomorphism, we can set $y_n=g_n$, and $\phi$ is entirely
determined by the formal series (the Nica-Speicher Fourier transform
\cite{NS})

\begin{equation}
\Phi(t) = \phi(g(t))=\sum_{n\ge 0}a_nt^{n}.
\end{equation}
Since $g(t)=\sum_{n\ge 0}t^nh_ng(t)^n$, we have
\begin{equation}
\Phi(t)
 = \sum_{n\ge 0}t^n\phi (h_n)\phi(g(t))^n
 = \sum_{n\ge O}\alpha_n t^n\Phi(t)^n,\ \text{with $\alpha_n=\phi(h_n)$.}
\end{equation}
Let $\psi$ be another multiplicative function such that $\psi(g_n)=b_n$ and
$\psi(h_n)=\beta_n$, and $\Psi(t)=\psi(g(t))$. Their convolution
$\eta=\phi\star\psi$ is determined by

\begin{align}
H(t)&= \phi\star\psi(g(t)) = (\phi\otimes\psi)\Delta g(t))\\
&= \left(\sum_{k\ge 0}\alpha_k t^kH(t)^k\right)
   \left(\sum_{\ell\ge 0}\beta_\ell t^lH(t)^\ell\right) \\
&=\sum_{n\ge 0}t^n\left(\sum_{k+\ell=n}\alpha_k\beta_\ell\right)H(t)^n.
\end{align}
Thus, convolution corresponds to the ordinary product of formal series
\begin{equation}
\hat\phi(t)=\sum_{n\ge 0}\alpha_nt^n,\ 
\hat\psi(t)=\sum_{n\ge 0}\beta_nt^n,\ 
\hat\eta(t)=\hat\phi(t)\hat\psi(t)=\sum_{n\ge 0}\gamma_n t^n
\end{equation}

since $H(t)$ satisfies
\begin{equation}
H(t) = \sum_{n\ge 0}\gamma_nt^nH(t)^n.
\end{equation}

As an illustration, the M\"obius function of the $\NC_{n+1}$ is the convolution
inverse of the $\zeta$ function, which is defined by $\zeta([\sigma,\pi])=1$
if $\sigma\le\pi$ and $0$ otherwise.
It is therefore characterized by $\zeta([0_{n+1},1_{n+1}])=\zeta(y_n)=1$.

If $\zeta(h_n)=\alpha_n$, then
\begin{equation}
Z(t) = \frac1{1-t}
     = \sum_{n\ge 0}\zeta(h_n)t^nZ(t)^n
     = \sum_{n\ge 0}\alpha_n\left(\frac{t}{1-t}\right)^n
\end{equation}
yields $\hat\zeta(t)=1+t$. Hence,
\begin{equation}
\hat\mu(t)=\frac1{1+t} \ \text{and}\ M(t)=\frac{1}{1+tM(t)},
\end{equation}
so that
\begin{equation}
M(t)=\frac{-1+\sqrt{1+4t}}{2t}.
\end{equation}

One can also count intervals and multichains.
Set $\zeta_k=\zeta^{\star k}$.
Then, $\hat\zeta_k(t)=(1+t)^k$. Hence $Z_k(t)$ satisfies
\begin{equation}
Z_k(t) = (1+tZ_k(t))^k,
\end{equation}
or alternatively 
\begin{equation}
X_k(t)
       = 1+tX_k(t)^k\ \text{with $X_k(t)=1+tZ_k(t)$},
\end{equation}
and we recover the fact that multichains of length $k$ are in bijection with
$(k+1)$-ary trees \cite{Ed82}.

In \cite{Ed80}, Edelman obtains a formula for the number of chains with
prescribed ranks $0_{n+1}<\pi_1<\ldots<\pi_r<\pi_{r+1}=1_{n+1}$. To derive it,
one can compute
\begin{equation}
\psi:=\varphi_{u_1}\star\varphi_{u_2}\star\cdots\star\varphi_{u_{r+1}},
\end{equation}
where $\varphi_u(g_n)=u^n$.
Then, $\hat\varphi_u(t)=1+tu$ and 
\begin{equation}
\hat\psi(t) = (1+tu_1)(1+tu_2)\cdots(1+tu_{r+1})=\lambda_t(U).
\end{equation}
Lagrange inversion yields
\begin{equation}
\psi(g_n)=\frac1{n+1}e_n[(n+1)U]
\end{equation}
and extracting the coefficient of a monomial, we obtain the number of chains
such that $\rk(\pi_i)-\rk(\pi_{i-1})=s_i$ is equal to
\begin{equation}
\frac1{n+1}\binom{n+1}{s_1}
\binom{n+1}{s_2}\cdots\binom{n+1}{s_{r+1}}.
\end{equation}
This is
\begin{equation}
\frac1{n+1}[u_1^{s_1}u_2^{s_2}\cdots u_{r+1}^{s_{r+1}}]\lambda_1[u_1+\cdots+u_{r+1}]^{n+1}
\end{equation}
which is equal to the coefficient of $m_\mu$ in $\omega(g)$, where $\mu$ is
the partition obtained by reordering the $s_i$, {\it i.e.}, to the scalar
product $\<e_\mu,g\>$. 

This last expression can be interpreted in terms of the Farahat-Higman
algebra. Let $c_\mu$ be the dual basis of $g^\mu$ in $Sym$ ({\it i.e.}
$c_\mu(-X)$ is what is denoted by $g_\mu$ in Macdonald's book \cite[Ex. 24-25,
p. 131-133]{Mcd}).
Then, the elementary symmetric functions are
\begin{equation}
e_k = \sum_{\kappa\vdash k}c_\kappa.
\end{equation}
Indeed, 
\begin{equation}
\<e_k,g^\kappa\>=\prod_i \frac{e_{\kappa_i}[\kappa_i+1]}{\kappa_i+1}=1
 \end{equation}
for all $\kappa\vdash k$.
Thus, it represents the sum of all permutations which can be written as a
minimal product of $k$ transpositions. Identifying $\NC_{n+1}$ with the
interval $[id_{n+1},(12\cdots n+1)]$ of the Cayley graph of $\SG_{n+1}$ as in
\cite{Bi2}, noncrossing partitions are identified with the permutations lying
on the minimal paths between the identity and the full cycle, the rank being
the transposition length.
If $\mu=(s_1,\ldots,s_{r+1})$, the scalar product $\<e_\mu,g_n\>$ is equal to
the coefficient of $c_n$ in the product $e_{s_1}s_{s_2}\cdots e_{s_{r+1}}$,
hence to the number of factorisations of the full cycle into a product of
permutations minimally factorisable into $s_1,s_2,\ldots$ transpositions, that
is, to the number of chains of noncrossing partitions with the prescribed
ranks.

As another example, since $c_n=M_n=p_n$, we can recover a result of Biane
\cite{Bi2}: the number of minimal factorizations of an $n$-cycle into a
product of cycles of orders $a_1,\ldots,a_r$ is the coefficient of $c_{n-1}$
in the product $c_{a_1-1}c_{a_2-1}\cdots c_{a_r-1}$, that is,
\begin{equation}
\<p_{a_1-1}\cdots p_{a_r-1},g_{n-1}(X)\>
= \<p_{a_1-1}\cdots p_{a_r-1},\frac1nh_{n-1}(nX)\>
= n^{r-1}.
\end{equation}

\section{Appendix}
\label{sec-appendix}

\subsection{Generating compatible pairs of compositions}

Given a composition $I$, the list of compositions $J$ compatible with $I$ 
can be computed as follows.

The composition whose descent set is the complement of the descent set of $I$
is ${\bar I}\tilde~$, the mirror conjugate of $I$.
Then, since we required that the sorted concatenation of the descent sets of
$I$ and $J$ form a word greater componentwise than the sorted concatenation of
$I$ and $\bar I^\sim$, the $J$ that are compatible with $I$ are those obtained
from ${\bar I}\tilde~$ by iterating the following process: given
$C=(c_1,\dots,c_n)$, for any $i>1$ such that $c_i>1$, change $C$ into $C'$ by
adding $1$ to $c_{i-1}$ and subtracting $1$ to $c_i$.

In particular, the set of compositions compatible with $I$ is equipped with a
natural order, its top element being ${\bar I}\tilde~$ and its bottom element
being $(|I|-k+1,1^{k-1})$ if $k$ is the length of ${\bar I}\tilde~$.

\medskip
{\footnotesize
For example, with $I=321$, its reverse conjugate is $1122$ and the whole list
of possibilities for $J$ contains $9$ elements which  can be drawn on the
following diagram
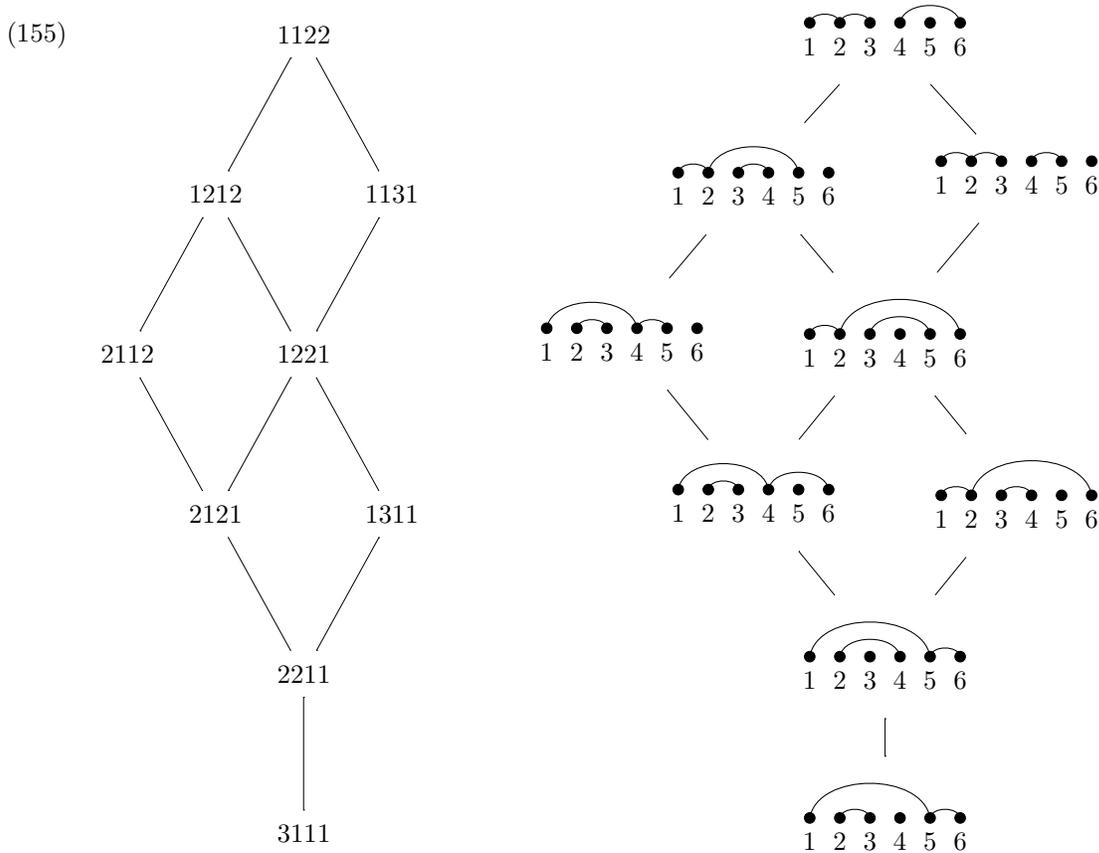
\begin{equation}
\setlength\unitlength{1.7mm}
\xymatrix@R=1.5cm@C=.7mm{
 & &
   *{\GrTeXBox{1122}}\arx1[dl]\arx1[dr] & \\
 & *{\GrTeXBox{1212}}\arx1[dl]\arx1[dr] & &
   *{\GrTeXBox{1131}}\arx1[dl] \\
   *{\GrTeXBox{2112}}\arx1[dr] & &
   *{\GrTeXBox{1221}}\arx1[dl]\arx1[dr] \\
 & *{\GrTeXBox{2121}}\arx1[dr] & &
   *{\GrTeXBox{1311}}\arx1[dl] \\
&& *{\GrTeXBox{2211}}\arx1[d] \\
&& *{\GrTeXBox{3111}} \\
}
\hskip1cm
\xymatrix@R=0.5cm@C=-12.7mm{
 & &
   *{\GrTeXBox{\Matchingsmall{6}{1/2,2/3,4/6}}}\arx1[dl]\arx1[dr] & \\
 & *{\GrTeXBox{\Matchingsmall{6}{1/2,2/5,3/4}}}\arx1[dl]\arx1[dr] & &
   *{\GrTeXBox{\Matchingsmall{6}{1/2,2/3,4/5}}}\arx1[dl] \\
   *{\GrTeXBox{\Matchingsmall{6}{1/4,4/5,2/3}}}\arx1[dr] & &
   *{\GrTeXBox{\Matchingsmall{6}{1/2,2/6,3/5}}}\arx1[dl]\arx1[dr] \\
 & *{\GrTeXBox{\Matchingsmall{6}{1/4,4/6,2/3}}}\arx1[dr] & &
   *{\GrTeXBox{\Matchingsmall{6}{1/2,2/6,3/4}}}\arx1[dl] \\
&& *{\GrTeXBox{\Matchingsmall{6}{1/5,5/6,2/4}}}\arx1[d] \\
&& *{\GrTeXBox{\Matchingsmall{6}{1/5,5/6,2/3}}} \\
}
\end{equation}
}

\subsection{Descent words and Motzkin paths}

We have seen that the compatible compositions are those whose concatenation
of descent sets are greater than $1\dots n$. The map sending compatible
compositions to such descent words is of course highly non injective and its
image set consists in the sorted words $s$ such that $s_i\geq i$, $s_i\leq n$,
and no value can be taken more than twice. Let us denote by $S_n$ this set of
words.

\medskip
{\footnotesize
For example, with $n=4$, we get the word $123$ eight times, and all other
words $133$, $223$ and $233$ twice each, for a total of $14$.
For general $n$, the number of pairs of compositions with a given word $s$ as
image is obviously $2^k$ where $k$ is the number of values used only once in
$s$. Moreover, if one counts the number of words by their number of doubled
letters (so that the first column is $1$ and represents $s=1\dots n$ with no
repeated letters), we find the following triangle:

\begin{equation}
\begin{array}{ccccccc}
  1 &       &      &     &     &     \\
  1 &       &      &     &     &     \\
  1 &     3 &      &     &     &     \\
  1 &     6 &    2 &     &     &     \\
  1 &    10 &   10 &     &     &     \\
  1 &    15 &   30 &   5 &     &     \\
  1 &    21 &   70 &  35 &     &     \\
  1 &    28 &  140 & 140 &  14 &     \\
\end{array}
\end{equation}
which is Sequence~A055151 of~\cite{Slo}, the triangular array of Motzkin paths
of length $n$ and with $k$ up steps.

In one wants to see how powers of $2$ come into play, one has to represent the
table as follows:
\begin{equation}
\begin{array}{cccccccccc}
  1 &       &      &     &     &     &     &     \\
    &     1 &      &     &     &     &     &     \\
  1 &       &    1 &     &     &     &     &     \\
    &     3 &      &   1 &     &     &     &     \\
  2 &       &    6 &     &   1 &     &     &     \\
    &    10 &      &  10 &     &   1 &     &     \\
  5 &       &   30 &     &  15 &     &  1  &     \\
    &    35 &      &  70 &     &  21 &     &  1  \\
 14 &       &  140 &     & 140 &     & 28  &  &1 \\
\end{array}
\end{equation}
Here, column $k$ corresponds to the number of words appearing $2^k$ times.
For example, the fifth line reads $2\cdot 2^0+6\cdot2^2+2^4=42$.
}

\medskip
We shall prove that $S_n$ is indeed equinumerous with Motzkin paths,
even with the extra parameter introduced above, but it will be easier to work
with the set $S'_n$ defined as the image of $S_n$ by the map 
\begin{equation}
w = w_1\dots w_n \mapsto (n+1-w_n)\dots (n+1-w_1).
\end{equation}
The condition on the words of $S_n$ translates in $S'_n$ as $w_i\leq i$, so
that the $w_i$ are parking functions. Now, the classical bijection between
nondecreasing parking functions and noncrossing partitions  sends a
noncrossing partition $c$ to the nondecreasing word where $i$ appears as many
times as the cardinality of the $i$-th part of $c$.

So $S'_n$ corresponds to the noncrossing partitions with parts at most $2$.
Read such a noncrossing partition from left to right and draw an up step if we
begin a part with two elements, a down step if we close such a part, and a
horizontal step if we have a singleton. This is the natural bijection between
these particular noncrossing partitions and Motzkin paths. Moreover, the
statistic of the number of repeated letters is sent to the number of parts
with two elements in the noncrossing partition, and then to the number of up
steps in the Motzkin path.

We then have
\begin{proposition}
The set $S'_n$ and Motzkin paths $M_n$ of $n$ are equinumerous and the
statistic of the number of repeated letters in $S'$ corresponds to the number
of up steps in $M_n$.
\end{proposition}

\medskip
{\footnotesize
For example, here is the whole of sublist $S_5$ consisting in words with two
pairs of repeated letters (10 elements) and their successive images by the
bijections.
}
\begin{equation}
\begin{split}
34455\ & \longleftrightarrow\ 11223\
        \longleftrightarrow\ \raisebox{-1em}{\Matchingsmall{5}{1/5,2/4}}
        \longleftrightarrow\ \
\begin{picture}(8,6) \up00\up22\hori44\down64\down82\pt{10}0 \end{picture}\\
24455\ & \longleftrightarrow\ 11224\
        \longleftrightarrow \raisebox{-1em}{\Matchingsmall{5}{1/5,2/3}}
        \longleftrightarrow\ \
\begin{picture}(8,6) \up00\up22\down44\hori62\down82\pt{10}0 \end{picture}\\
14455\ & \longleftrightarrow\ 11225\
        \longleftrightarrow\ \raisebox{-1em}{\Matchingsmall{5}{1/4,2/3}}
        \longleftrightarrow\  \
\begin{picture}(8,6) \up00\up22\down44\down62\hori80\pt{10}0 \end{picture}\\
33455\ & \longleftrightarrow\ 11233\
        \longleftrightarrow\ \raisebox{-1em}{\Matchingsmall{5}{1/5,3/4}}
        \longleftrightarrow\  \
\begin{picture}(8,6) \up00\hori22\up42\down64\down82\pt{10}0 \end{picture}\\
22455\ & \longleftrightarrow\ 11244\
        \longleftrightarrow\ \raisebox{-1em}{\Matchingsmall{5}{1/3,4/5}}
        \longleftrightarrow\ \
\begin{picture}(8,6) \up00\hori22\down42\up60\down82\pt{10}0 \end{picture}\\
23355\ & \longleftrightarrow\ 11334\
        \longleftrightarrow\ \raisebox{-1em}{\Matchingsmall{5}{1/2,3/5}}
        \longleftrightarrow\ \
\begin{picture}(8,6) \up00\down22\up40\hori62\down82\pt{10}0 \end{picture}\\
13355\ & \longleftrightarrow\ 11335\
        \longleftrightarrow\ \raisebox{-1em}{\Matchingsmall{5}{1/2,3/4}}
        \longleftrightarrow\ \
\begin{picture}(8,6) \up00\down22\up40\down62\hori80\pt{10}0 \end{picture}\\
22355\ & \longleftrightarrow\ 11344\
        \longleftrightarrow\ \raisebox{-1em}{\Matchingsmall{5}{1/2,4/5}}
        \longleftrightarrow\ \
\begin{picture}(8,6) \up00\down22\hori40\up60\down82\pt{10}0 \end{picture}\\
33445\ & \longleftrightarrow\ 12233\
        \longleftrightarrow\ \raisebox{-1em}{\Matchingsmall{5}{2/5,3/4}}
        \longleftrightarrow\ \
\begin{picture}(8,6) \hori00\up20\up42\down64\down82\pt{10}0 \end{picture}\\
22445\ & \longleftrightarrow\ 12244\
        \longleftrightarrow\ \raisebox{-1em}{\Matchingsmall{5}{2/3,4/5}}
        \longleftrightarrow\ \
\begin{picture}(8,6) \hori00\up20\down42\up60\down82\pt{10}0 \end{picture}\\
\end{split}
\end{equation}

Following~\cite{Slo}, there is a simple formula for $|S_{n,k}|$, the number of
elements of $S_n$ with $k$ repeated values: $|(S_{n,k})| = \binom{n}{2k}C_k$,
so that the cardinality of the set of compatible pairs of compositions is
\begin{equation}
\sum_{k\geq0} 2^{n-2k} \binom{n}{2k}C_k = C_{n+1},
\end{equation}
thanks to Touchard, cited by several authors on the Catalan webpage
of~\cite{Slo}.
So we have proved by a simple and meaningful bijection that indeed $\Delta
G_n$ has Catalan terms.

\bigskip
{\bf Acknowlegements. } This research has been partially supported by the
project CARPLO of the Agence Nationale de la recherche (ANR-20-CE40-0007).



\footnotesize
\begin{thebibliography}{aa}
%
\bibitem{Bi1}{ P. Biane},
{ Some properties of crossings and partitions},
Discrete Math. { 175} (1997), 41--53.
%
\bibitem{Bi2}{ P. Biane},
{ Minimal factorizations of a cycle and central multiplicative functions on
the infinite symmetric group},
JCTA { 76} (1996), 197--212.
%
\bibitem{BFK}{ C. Brouder, A. Frabetti, C. Krattenthaler},
{ Non-commutative Hopf algebra of formal diffeomorphisms},
Adv. Math.   200  (2006),  479--524.
%
%
\bibitem{Ed80}{ P. H. Edelman}, 
{ Chain enumeration and noncrossing partitions},
Discrete Math. { 31} (1980), 171--180.
%
\bibitem{Ed82}{ P. H. Edelman}, 
{ Multichains, noncrossing partitions and trees}, 
Discrete Math. { 40} (1982), 171--179.
%
\bibitem{EH}{ R. Ehrenborg, A.Happ},
{ The antipode of the noncrossing partition lattice},
Advances in Applied Mathematics,
110 (2019),  76--85.
%
\bibitem{Einz}{ H. Einziger},
{ Incidence Hopf Algebras: Antipodes, Forest Formulas, and Noncrossing
Partitions},
Ph.D. Thesis, The George Washington University, 2010.
%
\bibitem{NCSF1}{ I. M. Gelfand, D. Krob, A. Lascoux, B. Leclerc,
V.~S. Retakh,  J.-Y. Thibon},
{ Noncommutative symmetric functions},
Adv. in Math. 112 (1995), 218--348.
%
%
\bibitem{GJ1}{ I. P. Goulden, D. M. Jackson},
{ The Combinatorial Relationship Between Trees, Cacti and Certain
Connection Coefficients for the Symmetric Group},
Europ. J. Combinatorics { 13} (1992), 357--365.
%
\bibitem{GJ2}{ I. P. Goulden, D. M. Jackson},
{ Symmetric functions and Macdonald's result for top connexion coefficients in
the symmetric group},
J. Algebra { 166} (1994), no. 2, 364--378.
%
\bibitem{Hai1}{ M. Haiman},
{ Conjectures on the quotient ring by diagonal invariants},
J. Algebraic Combin.  3 (1994), 17--36.
%
\bibitem{HNTtrees}{F. Hivert, J.-C. Novelli, J.-Y. Thibon},
Trees, functional equations, and combinatorial Hopf algebras,
European J. Combin. 29 (2008), no. 7, 1682--1695.
%
\bibitem{JMNT}{ M. Josuat-Verg\`es, F. Menous, J.-C. Novelli, J.-Y. Thibon},
{ Free cumulants, Schröder trees, and operads},
Advances in Applied Mathematics { 88}  (2017),  92--119.
%
%
\bibitem{Kr}{ G.  Kreweras},
{ Sur les partitions non croisées d'un cycle (French)},
Discrete Math. 1 (1972), 333--350.
%
\bibitem{Mcd}{ I. G. Macdonald},
{ Symmetric functions and Hall polynomials},
2nd ed., Oxford University Press, 1995.
%
\bibitem{NS}{ A. Nica, R. Speicher},
{ A ``Fourier transform'' for multiplicative functions on noncrossing
partitions},
J. Algebraic Combin. 6 (1997),  141--160
%
\bibitem{NTLag}{ J.-C Novelli}, { J.-Y. Thibon},
{ Noncommutative symmetric functions and Lagrange inversion},
Adv. Appl. Math. { 40} (2008), 8--35.
%
\bibitem{NTDup}{ J.-C. Novelli}, { J.-Y. Thibon},
{ Duplicial algebras and Lagrange inversion}, arXiv:1209.5959.
%
\bibitem{NTm}{ J.-C. Novelli}, { J.-Y. Thibon},
{ Hopf Algebras of m-permutations, (m+1)-ary trees, and m-parking functions},
Advances in Applied Mathematics { 117} (2020) 102019.
%
\bibitem{NTpark}{ J.-C. Novelli}, { J.-Y. Thibon},
{ Hopf algebras and dendriform structures arising from parking functions},
Fund. Math.   193  (2007),   189--241.
%
\bibitem{Ran}{ G. N. Raney}, 
{ Functional composition patterns and power series reversion},
Trans. Amer. Math. Soc.  94 (1960), 441--451.
%
\bibitem{SU}{ R. Simion, D. Ullman},
{ On the structure of the lattice of noncrossing partitions},
Disc. Math.  98 (1991), 193--206.
%
\bibitem{Slo}{ N. J. A. Sloane},
{ The On-Line Encyclopedia of Integer Sequences},\\
http://www.research.att.com/~njas/sequences/.
%
\bibitem{Spei1}{R. Speicher},
{ Multiplicative functions on the lattice of noncrossing partitions and free
convolution},
Math. Ann. 298  (1994),  611--628.
%
%
%
\end{thebibliography}
\end{document}